\documentclass[12pt]{article}
\usepackage{graphicx}
\usepackage{amsmath,amsthm,amssymb,enumerate}
\usepackage{euscript,mathrsfs}
\usepackage[left=1.5cm,right=1.5cm,top=3.5cm,bottom=3.5cm]{geometry}
\usepackage{url}
\usepackage{color}
\catcode`\@=11 \@addtoreset{equation}{section}

\catcode`\@=12
\allowdisplaybreaks

\newtheorem{Theorem}{Theorem}[section]
\newtheorem{Proposition}[Theorem]{Proposition}
\newtheorem{Lemma}[Theorem]{Lemma}
\newtheorem{Corollary}[Theorem]{Corollary}

\theoremstyle{definition}
\newtheorem{Definition}[Theorem]{Definition}

\newtheorem{Remark}[Theorem]{Remark}

\newcommand{\bTheorem}[1]{
\begin{Theorem} \label{T#1} }
\newcommand{\eT}{\end{Theorem}}

\newcommand{\bProposition}[1]{
\begin{Proposition} \label{P#1}}
\newcommand{\eP}{\end{Proposition}}

\newcommand{\bLemma}[1]{
\begin{Lemma} \label{L#1} }
\newcommand{\eL}{\end{Lemma}}

\newcommand{\bCorollary}[1]{
\begin{Corollary} \label{C#1} }
\newcommand{\eC}{\end{Corollary}}

\newcommand{\bRemark}[1]{
\begin{Remark} \label{R#1} }
\newcommand{\eR}{\end{Remark}}

\newcommand{\bDefinition}[1]{
\begin{Definition} \label{D#1} }
\newcommand{\eD}{\end{Definition}}

\newcommand{\Del}{\Delta_x}

\DeclareMathOperator{\dist}{dist}

\newcommand{\vrd}{\varrho_{\delta}}
\newcommand{\vud}{\vu_{\delta}}

\newcommand{\bfphi}{\boldsymbol{\varphi}}

\newcommand{\bFormula}[1]{
\begin{equation} \label{#1}}
\newcommand{\eF}{\end{equation}}

\newcommand{\vr}{\varrho}
\newcommand{\vre}{\vr_\ep}

\newcommand{\vue}{\vu_\ep}

\newcommand{\vu}{\vc{u}}
\newcommand{\vc}[1]{{\bf #1}}

\newcommand{\Div}{{\rm div}_x}
\newcommand{\Grad}{\nabla_x}

\newcommand{\tn}[1]{\mathbb{#1}}
\newcommand{\dx}{\,{\rm d} {x}}

\newcommand{\intO}[1]{\int_{\Omega} #1 \ \dx}

\newcommand{\ep}{\varepsilon}

\newcommand{\R}{\mathbb{R}}

\definecolor{Cgrey}{rgb}{0.85,0.85,0.85}
\definecolor{Cblue}{rgb}{0.50,0.85,0.85}
\definecolor{Cred}{rgb}{1,0,0}
\definecolor{fancy}{rgb}{0.10,0.85,0.10}

\newcommand\Cbox[2]{%
    \newbox\contentbox%
    \newbox\bkgdbox%
    \setbox\contentbox\hbox to \hsize{%
        \vtop{
            \kern\columnsep
            \hbox to \hsize{%
                \kern\columnsep%
                \advance\hsize by -2\columnsep%
                \setlength{\textwidth}{\hsize}%
                \vbox{
                    \parskip=\baselineskip
                    \parindent=0bp
                    #2
                }%
                \kern\columnsep%
            }%
            \kern\columnsep%
        }%
    }%
    \setbox\bkgdbox\vbox{
        \color{#1}
        \hrule width  \wd\contentbox %
               height \ht\contentbox %
               depth  \dp\contentbox
        \color{black}
    }%
    \wd\bkgdbox=0bp%
    \vbox{\hbox to \hsize{\box\bkgdbox\box\contentbox}}%
    \vskip\baselineskip%
}

\date{}

\begin{document}


\title{Dissipative solutions to compressible Navier-Stokes equations with general inflow-outflow data:
existence, stability and weak strong uniqueness}

\author{Young-Sam Kwon
\thanks{The work of Young-Sam Kwon was supported by Basic Science Research Program through the National Research Foundation of
Korea(NRF) funded by the Ministry of
Education(NRF-
2017R1D1A1B03030249)   }
\and Antonin Novotny
 \and Vladyslav Satko
}

\date{\today}

\maketitle

\bigskip

\centerline{Department of Mathematics,  Dong-A University,}

\centerline{Busan 604-714, Republic of Korea}


\centerline{and}

\centerline{IMATH, EA 2134, Universit{\' e} du Sud Toulon-Var}

\centerline{BP 20132, 83957 La Garde, France}

\centerline{and}

\centerline{IMATH, EA 2134, Universit{\' e} du Sud Toulon-Var}

\centerline{BP 20132, 83957 La Garde, France}

\bigskip

\begin{abstract}
So far existence of dissipative weak solutions for the compressible Navier-Stokes equations (i.e. weak solutions satisfying the relative energy inequality) is known only
in the case of boundary conditions with non zero inflow/outflow (i.e., in particular, when the normal component of the velocity on the boundary of the flow domain is equal to zero).
Most of physical applications (as flows in wind tunnels, pipes, reactors of jet engines) requires to consider non-zero inflow-outflow boundary condtions.

We prove existence of dissipative weak solutions to the compressible Navier-Stokes equations in barotropic regime (adiabatic coefficient $\gamma>3/2$, in three dimensions, $\gamma>1$ in two dimensions) with large velocity prescribed at the boundary and
large density prescribed at the inflow boundary of a bounded piecewise regular Lipschitz domain, without any restriction neither on the shape  of the inflow/outflow boundaries nor on the shape  of the domain.

It is well known that the relative energy inequality has many applications, e.g., to investigation of incompressible or inviscid limits, to the dimension reduction of flows, to the error estimates of numerical schemes. In this paper we deal with one of its basic applications, namely weak-strong uniqueness principle.

\end{abstract}



\bigskip

\section{Introduction}
\label{i}

 We consider the system of equations governing the
non steady motion of a compressible viscous fluid driven by general in/out flow boundary conditions on general bounded domains.
The mass density $\vr = \vr(t,x)$ and
the velocity $\vu = \vu(t,x)$, $(t,x)\in I\times\Omega\equiv Q_T$, $I=(0,T)$ of the fluid satisfy the Navier--Stokes system,
\begin{eqnarray}
\label{i1}
\partial_t\vr+\Div (\vr \vu) &=& 0,\\
\label{i2}
\partial_t(\vr\vu)+\Div (\vr \vu \otimes \vu) + \Grad p(\vr) &=& \Div \mathbb{S}(\Grad \vu),
\end{eqnarray}
in $\Omega \subset R^d$, $d=2,3$, where  the stress tensor is defined by
 \[
 \mathbb{S}(\Grad \vu) = \mu \left( \Grad \vu + \Grad^t \vu \right) + \lambda \Div \vu \mathbb{I}, \ \mu > 0 , \ \lambda \geq 0,
 \]
 and $p = p(\vr)$ is the barotropic pressure.

The system is completed with initial conditions
\begin{equation}\label{initc}
\vr(0)=\vr_0,\quad\vr\vu(0)=\vr_0\vu_0
\end{equation}
and boundary conditions
\begin{equation}
\label{i4}
\vc{u} |_{\partial \Omega} = \vu_B, \ \vr|_{\Gamma_{\rm in}} = \vr_B,
\end{equation}
where
\begin{equation}
\label{i5}
\Gamma_{\rm in} = \left\{ x \in \partial \Omega \ \Big| \ {\vc{u}_B} \cdot \vc{n} < 0 \right\},\
\Gamma_{\rm out} = \left\{ x \in \partial \Omega \ \Big| \ {\vc{u}_B} \cdot \vc{n} > 0 \right\}.
\end{equation}
In the above $\vc n$ is the outer normal to the boundary $\partial\Omega.$


The investigation and the better insight to the equations in this setting   is important for many real world applications. In fact, this is a natural and
basic abstract setting
for flows in pipelines, wind tunnels, turbines to name a few concrete examples.
{Numerical modeling of fluid flow in the portions of cooling circuits of nuclear power stations, in the gas transporting industrial pipelines,
in the reactors of jet engines and in many other situations, requires this or similar boundary value setting rather than academic no-slip, Navier or periodic
boundary conditions.}

 Existence of (renormalized bounded energy) weak solutions for system (\ref{i1}--\ref{i5}) is well known  -for the pressure $p(\vr)$ behaving as $\vr^\gamma$ at infinity with $\gamma>d/2$- in the "simple" case of zero inflow and outflow boundary data with no-slip
or slip (Navier) boundary conditions (when, in particular, $\vc u\cdot\vc n|_{\partial\Omega}=0$) since the end of the last century/beginning of this
century, cf. \cite{FNP}, Feireisl\cite{FEnonmon} and monographs by Lions \cite{LI4}, Feireisl \cite{EF70} and \cite{NOST} (see also
an alternative approach by Bresch, Jabin \cite{BreJab} working for $\gamma\ge 9/5$ ($d=3$) with possibly non monotone pressure). In this situation, it is also known, that any bounded energy weak solution is dissipative, meaning that it obeys the so called relative energy inequality and consequently satisfies, in particular, the weak strong uniqueness principle, see \cite{FeJiNo}, \cite{FNSun}  (and the seminal paper of Dafermos \cite{Dafermos} for the general introduction of relative energy (entropy)
method in the fluid mechanics).

Existence of (renormalized bounded energy) weak solutions for system (\ref{i1}--\ref{i5}) with large boundary data exhibits many additional difficulties.
It was investigated recently in Novo \cite{NOVO}, Girinon \cite{Girinon} with several geometrical restrictions and in \cite{JiNoSIMA}, \cite{ChoeNoY} in full generality.
If the initial data and the boundary are smooth, the same problem admits local in time strong solutions which become global if the initial data are sufficiently small,
see Valli, Zajaczkowski \cite{VAZA}. A general question arises whether strong solutions are unique in the class of weak solutions, at least on the lifespan of the former.
In contrast to the situation with no inflow/outflow, in this general setting, it is not known, whether the class of renormalized bounded energy weak solution coincides with the class of dissipative solutions, i.e., whether the weak strong uniqueness is true. This question remains apparently an interesting open problem.

We are, however, able to construct a subclass of weak solutions, called {\it dissipative} {\it weak solutions},  that obey, in particular, the weak strong uniqueness principle. The theorem (and proof) of the existence of (renormalized bounded energy)  dissipative weak solutions is the first main goal of the present paper.

Except its important application potential, there are notably two features that distinguish our result from the similar result known for the no-slip (or Navier) boundary conditions:
\begin{enumerate}
\item The test density in the relative energy inequality cannot be taken arbitrarily but must obey the continuity equation with transporting velocity which is an arbitrary vector field satisfying the boundary conditions of the problem, and which is at the same time the test velocity in the relative energy inequality.
\item Presuure $p(\vr)$ is only $L^1$-integrable near the boundary. In absence of higher integrability, one must show that it is "equi-integrable" near the boundary.
This result is formulated in Lemma \ref{localp} and it is of independent interest.
\end{enumerate}

As in the case of zero boundary conditions, this result (in particular, the relative energy inequality) opens the way to many
applications in geometrical setting with non zero inflow/outflow that are really interesting from the engineering and physical points of view both
on the theoretical level
(as e.g. rigorous investigation of model reduction through singular distinguished limits that plays crucial role in derivation of simplified models in the physics of atmosphere, cf. Klein et al. \cite{Kleinetal}, in the spirit of \cite[2nd. edition, Chapter 9]{FEINOV}) and on the level of numerical analysis (as e.g. derivation of rigorous unconditional error estimates in the spirit of \cite{GaHeMaNo}, \cite{GaMaNo}), which were, so far, out of reach of the analysis.

In this paper, we develop one of them, namely, we prove stability of strong solutions in the class of dissipative weak solutions and the weak-strong uniqueness principle for the dissipative weak solutions. This is the second main result of the present paper.

The paper is organized as follows. In Section \ref{M} we expose the definition of dissipative weak solutions and state the main theorems
(Theorem \ref{TM1} about the existence of weak dissipative solutions and Theorem \ref{TWS1} about the stability of strong solutions in the
class of weak solutions and the weak strong uniqueness). The construction of solutions is explained in Section \ref{A} which provides in Lemma
\ref{LLemmaA} existence of generalized solutions to the approximate system. This lemma is proved in Section \ref{4}. This is the real starting point
of the topic. A relative energy inequality for the approximate problem is derived in Section \ref{5}, see Lemma \ref{LRenergy}. It is the starting
point of the proof of Theorem \ref{TM1} which is performed in Sections \ref{D} and \ref{sE}. Section \ref{WS} is devoted to the proof of
Theorem \ref{TWS1}. Finally, in Section \ref{CR} we show that the results of Theorems \ref{TM1} and \ref{TWS1} can be extended to the piecewise regular Lipschitz domains and to certain nonmonotone pressure laws.

Throughout the paper, we use the standard notation for Sobolev and Bochner spaces, see, e.g., the book of Evans \cite{Evans}.

\section{Main results}
\label{M}

 In order to avoid additional technicalities, we suppose that the domain and boundary  data satisfy
\begin{equation} \label{M1}
\Omega\;\mbox{a bounded domain},\;\partial\Omega\in C^{2},\;
\vu_B \in C^2(\partial \Omega; R^d), \; 0< \underline\vr_B\le \vr_B\le \overline\vr_B,\;\vr_B \in C(\overline\Gamma_{\rm in}).
\end{equation}
The results however hold also for piecewise $C^2$ Lipschitz domains. This extension will be discussed in the last section.

Throughout the paper, we consider strictly increasing pressure satisfying at least
\bFormula{pressure}
{p}\in C[0,\infty)\cap C^1(0,\infty),\; {p}(0)= 0, \; p'(\vr)>0.
\eF
All results of this paper can be however extended to certain non-monotone pressure laws
(with non-monotonicity at least at a compact portion of the interval $[0,\infty)$).
We shall comment about this issue in the last section.

For further convenience, it will be useful to introduce the Helmholtz function
\bFormula{H}
H(\vr)=\vr\int_1^\vr\frac{p(z)}{z^2}{\rm d} z
\eF
and the relative energy function
\bFormula{E}
E(\vr|r)=H(\vr)-H'(r)(\vr-r)-H(r).
\eF
One can verify by the direct calculation that
\bFormula{H+}
\vr H'(\vr)-H(\vr)=p(\vr)\;\mbox{and consequently}\; H''(\vr)=\frac{p'(\vr)}\vr.
\eF

We begin with the definition of dissipative weak solutions to system (\ref{i1}--\ref{i5}). In this definition,
$\Omega$ is a bounded Lipschitz domain, $p\in  C[0,\infty)\cap C^1(0,\infty)$ and $\vu_B\in W^{1,\infty}(\partial\Omega)$,
$\vr_B\in L^\infty(\partial\Omega)$.
\\ \\
{\bf Definition \ref{M}.1} [Dissipative weak solutions to system (\ref{i1}--\ref{i5})]
\\
{\it We say that $(\vr, \vu)$ is a dissipative weak solution of problem (\ref{i1}--\ref{i5}) if:
\begin{description}
\item{1.} There exists a Lipschitz extension $\vu_\infty\in W^{1,\infty}(\Omega;\R^d)$ of $\vu_B$ whose divergence is non negative
in a certain interior neighborhood of $\partial\Omega$,  i.e.
\bFormula{hatU}
{\rm div}\vu_\infty\ge 0 \;\mbox{a.e. in}\; \hat{U}^-_h\equiv\{ x\in\Omega\,|\,{\rm dist}(x,\partial\Omega)<h\},\; h>0
\eF
and
\begin{align}\label{fs}
&\vr \in L^\infty (0,T; L^\gamma(\Omega))\;\mbox{with some $\gamma>1$},\; \ 0 \leq \vr \ \mbox{a.e. in}\ (0,T)\times \Omega,\nonumber\\
&
 p(\vr)\in L^1(Q_T), \;
\vc v:=\vu-\vu_\infty \in L^2(0,T; W_0^{1,2}(\Omega; R^d)).
\end{align}
\item{2.}
Function $\vr\in C_{\rm weak}([0,T], L^\gamma(\Omega))$\footnote{We say that $f\in  C_{\rm weak}([0,T], L^p(\Omega))$ iff $f:[0,T]\to L^p(\Omega)$ is defined everywhere on $[0,T]$, $f\in L^\infty(0,T; L^p(\Omega))$ and the map
$t\mapsto\intO{f\varphi(t,\cdot)}\in C[0,T]$ for
all $\varphi\in L^{p'}(\Omega)$.}   and the integral identity
\begin{align}\label{ce}
&\intO{\vr(\tau,\cdot)\varphi(\tau,\cdot)} - \intO{\vr_0(\cdot)\varphi(0,\cdot)}\nonumber\\
&
=\int_0^\tau\intO{\Big(\vr\partial_t\varphi + \vr \vu \cdot \Grad \varphi\Big) }{\rm d}t - \int_0^\tau\int_{\Gamma_{\rm in}} \vr_B \vu_B \cdot \vc{n} \varphi
\ {\rm d}S_x{\rm d} t
\end{align}
holds for any $\tau\in [0,T]$ and $\varphi \in C_c^1([0,T]\times({\Omega}\cup\Gamma_{\rm in}))$.
\item{3.} Function
$\vr\vu\in  C_{\rm weak}([0,T], L^{\frac{2\gamma}{\gamma +1}}(\Omega;R^d))$, and the integral identity
\begin{align}\label{me}
&\intO{\vr\vc v(\tau,\cdot)\cdot\bfphi(\tau,\cdot)} - \intO{\vr_0\vc v_0(\cdot)\bfphi(0,\cdot)}\nonumber\\
&=
\int_0^\tau\intO{  \Big(\vr\vc v\cdot\partial_t\varphi+ \vr \vu\cdot \Grad \bfphi \cdot \vu  -\vr\vu\cdot\Grad(\vu_\infty\cdot \bfphi) + p(\vr) \Div \bfphi
-  \mathbb{S}(\Grad \vu) : \Grad \bfphi\Big) }{\rm d}t
\end{align}
holds for any  $\tau\in [0,T]$ and any $\bfphi \in C^1_c ([0,T]\times\Omega; R^d)$.
\item{4.} { Function $\vr|\vu|^2\in L^\infty(I;L^1(\Omega))$ and the so called relative energy inequality}
\begin{align}\label{rei}
&\intO{\Big(\frac 12\vr|\vc v-\vc V|^2+E(\vr|r)\Big)(\tau)}+\int_0^\tau\intO{\tn S(\Grad\vu):\Grad(\vc v-\vc V)}{\rm d}t\nonumber\\
&
\le
{ \intO{\Big(\frac 12\vr_0|\vu_0-\vc U(0,\cdot)|^2+E(\vr_0|r(0,\cdot))\Big)}}\nonumber\\
&
\qquad+\int_0^\tau\int_{\Gamma_{\rm in}}\Big(H(r)-rH'(r)-H(\vr_B) +\vr_B H'(r)\Big)\vu_B\cdot\vc n{\rm d}S{\rm d}t
\nonumber\\
&\qquad+\int_0^\tau\intO{\Big(\vr(\vc V-\vc v)\cdot\partial_t\vc V +\vr\vu\cdot\Grad\vc U\cdot(\vc V-\vc v)\Big)}{\rm d}t
\nonumber\\
&\qquad+\int_0^\tau\intO{\Big( p(r)
 -p'(r)(r-\vr) -p(\vr)\Big){\rm div}\vc U }{\rm d}t
\nonumber\\
&\qquad+\int_0^\tau\intO{\Big(-p'(r)\vc v\cdot\Grad r- p(r){\rm div}\vc V+\frac{r-\vr }r p'(r)(\vc v-\vc V)\cdot\Grad r\Big)}{\rm d}t
\end{align}
holds with a.e. $\tau\in I$ and with any test functions $(r,\vc U)$,
\begin{equation}\label{Ugen}
\vc U= C^1(\overline Q),\; \vc U|_{\overline I\times\partial\Omega}=\vu_B, \;\vc V=:\vc U-\vu_\infty
\end{equation}
and
\begin{equation}\label{rgen}
0<r\in C^1(\overline{ Q_T})\;\,\mbox{satisfying}\;\;\partial_t r+{\rm div}(r\vc U)=0.
\end{equation}
\end{description}
}
{
\bRemark{R1-}

\begin{enumerate}
\item
A Lipschitz extension $\vu_\infty$ of $\vu_B$ verifying (\ref{hatU}) always exists. Indeed, according to \cite[Lemma 3.3]{Girinon}), for any
$\vc V\in W^{1,\infty}(\partial\Omega;\R^d)$ (where $\Omega\subset R^d$ is a bounded Lipschitz domain) there is $h>0$ and a vector field
\bFormula{Le}
\vc V_\infty\in W^{1,\infty}(\Omega),\quad {\rm div}{\vc V}_\infty\ge 0\;\mbox{ a.e. in $\hat{U}^-_h$}
\eF
verifying ${\vc V}_\infty|_{\partial\Omega}=\vc V$.

\item If $\Omega$ is bounded domain of class $C^2$  and $\int_{\partial\Omega}\vu_B\cdot\vc n {\rm d}S=0$ then $\vu_\infty$ can be chosen in such
a way that
$$
{\rm div}\vu_\infty=0\;\mbox{in $\Omega$},
$$
see Galdi \cite[Theorem IV.6.1]{Galdi}.
\item If $\partial\Omega$ is in class $C^2$, given $\vc U$ in class (\ref{Ugen}), the Cauchy-Lipschitz theory guarantees that equation  (\ref{rgen}) admits infinitely many solutions $r$ in class (\ref{rgen}) according to the free choice of the initial condition.
They can be constructed by the method of characteristics e.g. as follows: One extends $\vc U$ to $\tilde{\vc U}\in C^1_c(R\times R^d)$ and
finds its flow $X(\cdot;\cdot)$ (by definition $X(\cdot,x)$ is the (unique) solution of the ODE $y'(t)=\tilde{\vc U}(t,y)$,
$y(0)=x$). We know that $X\in C^1(R\times R^d)$, for all $t\in R$, $X(t;\cdot)$ is a bijection from $R^d$ to $R^d$ and $X(t+s,x)=
X(t;X(s;x))$, cf. \cite[Theorem 5.13]{BG10}. Consequently, for any given ${\mathfrak{R}}_0\in C^1(R^d)$ function ${\mathfrak{R}}(t,x)=
\mathfrak{R}_0(X(-t;x)){\rm exp}(-
\int_0^t{\rm div}\tilde{\vc U}(s,X(s-t;x)){\rm d}s$ solves the continuity equation $\partial_t\mathfrak{R}+{\rm div}(\mathfrak{R}\tilde{\vc U})=0$
in $R\times R^d$. Now, it is enough to take $r=\mathfrak{R}|_{[0,T]\times\overline\Omega}$.
\item  Equation (\ref{ce}) implies the total mass inequality
\bFormula{mi}
\intO{\vr(\tau)}\le\intO{\vr_0}-\int_0^\tau\int_{\Gamma_{\rm in}}\vr_B\vu_B\cdot\vc n{\rm d} S_x{\rm d}t
\eF
for all $\tau\in [0,T]$. To see it, it is enough to take for the test functions in (\ref{ce}) a convenient sequence $\varphi=\varphi_\ep$, $\ep>0$
as, e.g.,
\bFormula{L9}
\varphi_\ep(x) = \left\{ \begin{array}{l} 1 \ \mbox{if} \ {\rm dist}(x, \Gamma_{\rm out}) > \ep\\
\frac{1}{\ep} {\rm dist}(x, \Gamma_{\rm out}) \ \mbox{if}\ {\rm dist}(x, \Gamma_{\rm out}) \leq \ep  \end{array} \right\}
\eF
and send $\ep\to 0$.
\item Regularity of the test functions $r,\vc U$ in (\ref{Ugen}), (\ref{rgen}),  can be weaken, in particular, up to $r$, $\vc U$ continuous functions on
$\overline{Q_T}$ and $\partial_t(r,\vc U)$, $\Grad(r,\vc U)$ $\in L^2(I; C(\overline\Omega))$.
\item One can define {\it bounded energy weak solutions} requiring satisfaction of Items 1.-3. of Definition \ref{M1}, and energy inequality which reads,
\bFormula{ei}
\intO{\Big(\frac 12\vr|\vc v|^2+H(\vr)\Big)(\tau)}
+\int_0^\tau\intO{\tn S(\Grad\vc v):\Grad\vc v}{\rm d}t
\eF
$$
\le  \intO{\Big(\frac 12\vr_0|\vc v_0|^2+H(\vr_0)\Big)}
-\int_0^\tau\int_{\Gamma_{\rm in}} H(\vr_B)\vu_B\cdot\vc n{\rm d}S_x{\rm d}t
$$
$$
+
\int_0^\tau\intO{\Big(-p(\vr){\rm div}\vu_\infty-\tn S(\Grad\vu_\infty):\Grad\vc v-\vr\vu\cdot\Grad\vu_\infty\cdot\vc v
)}{\rm d}t,
$$
see \cite[Definition 2.1]{JiNoSIMA}.

We shall construct the dissipative weak solutions in such a way that they are also the bounded energy weak solutions, i.e., they satisfy, in addition to all items
in Definition \ref{M1}, also the energy inequality (\ref{ei}).

A natural question arises whether any dissipative weak solution is a bounded energy weak solution, i.e. if
 it satisfies also (\ref{ei}).  We recall that in the case of zero inflow/outflow boundary conditions,
the answer is "yes" and its proof takes one line, see \cite{FeJiNo}. In the case of non zero inflow/outflow boundary conditions
the answer is "yes", at least provided the function $p(z)/z^2$ is integrable near $0$. In this case, one can suppose without the loss of generality that $H'$ is a positive function. The process of deducing from the relative energy inequality (\ref{rei}) the energy inequality
(\ref{ei}) goes as follows: We take 1) in (\ref{ce}) the test function $\varphi=-\frac12{|\vc V|}^2$, 2) in (\ref{me}) the test function
$\varphi=\vc V$ and finally 3) in (\ref{ce}) the test function $\varphi=H'(r)\varphi_\ep$ (cf. (\ref{L9})). Adding the results of all three
above steps to (\ref{rei}), we obtain, after a long calculation (which uses identities (\ref{rgen}) and (\ref{H+})) and after
sending $\ep \to 0$, the energy inequality (\ref{ei}). The crucial point in this process is the treatment in the limit of terms
$\int_0^\tau\intO{H'(r)(\vu-\vu_\infty)\cdot \Grad\varphi_\ep}{\rm d}t$ and $\int_0^\tau\intO{H'(r)\vu_\infty)\cdot \Grad\varphi_\ep}{\rm d}t$ in step 3). 
In fact, $\lim_{\ep\to 0}$ of the first one is equal to $0$ by virtue of the Hardy inequality and $\limsup_{\ep\to 0}$ of the second one is non positive provided $H'(r)\ge 0$ (since $\Grad\dist(x, \Gamma_{\rm out})\to -\vc n(x_0)$ as $x\to x_0\in \Gamma_{\rm out}$).
\end{enumerate}
\eR
}
\noindent
{\bf Definition \ref{M}.2}
{\it We say that the couple $(\vr,\vu)\in C_{\rm weak}(\overline I;L^p(\Omega))\times L^2(0,T; W^{1,2}(\Omega,\R^d))$, $p>1$ is a
renormalized solution of the continuity equation
if $b(\vr)\in C_{\rm weak}([0,T];L^1(\Omega))$ (not relabeled in time) and if it satisfies in addition to the continuity equation (\ref{ce}) also equation
\begin{align} \label{P3}
&\intO{ (b(\vr) \varphi) (\tau)} -\intO{ b(\vr_0) \varphi(0)}\nonumber\\
&
=\int_0^\tau\intO{ \Big( b(\vr)\partial_t\varphi b(\vr) \vu \cdot \Grad \varphi -\varphi\left( b'(\vr) \vr - b(\vr) \right) \Div \vu \Big)}  {\rm d} t
- \int_0^\tau\int_{\Gamma_{\rm in}} b(\vr_B)
\vu_B \cdot \vc{n} \varphi \ {\rm d}S_x {\rm d} t
\end{align}
for any $\varphi \in C_c^1([0,T]\times({\Omega}\cup\Gamma_{\rm in}))$, and any continuously differentiable $b$ with $b'$ having a compact support
in $[0,\infty)$.

A (dissipative) weak solution to problem (\ref{i1}--\ref{i5}) satisfying in addition  renormalized continuity equation (\ref{P3})
is called a {\it renormalized (dissipative) weak solution}.
}


Our first main result is the following theorem.

\begin{Theorem}{\rm [Existence of dissipative weak solutions]} \label{TM1}
Let $\Omega \subset R^d$, $d = 2,3$ be a bounded domain of class { $C^{2}$.} Let the boundary data $\vu_B$, $\vr_B$ satisfy (\ref{M1}).
 Assume that the pressure satisfies hypotheses (\ref{pressure}) and
\bFormula{pressure1}
p'(\vr)\ge a_1\vr^{\gamma-1}-b,\,p(\vr)\le a_2\vr^\gamma+b,\;\gamma>d/2,\; a_1,a_2>0,\, b\ge 0.
\eF
Suppose finally that the initial data have the finite energy and the finite mass,
\bFormula{feid}
\intO{\Big(\frac 12\vr_0|\vu_0|^2+H(\vr_0)\Big)}<\infty,\quad 0\le\vr_0, \;\intO{\vr_0}>0.
\eF
{ Then for any Lipschitz extension $\vu_\infty$ of $\vu_B$ verifying (\ref{hatU})
problem (\ref{i1}--\ref{i5}) possesses at least one renormalized dissipative  weak solution $(\vr, \vu)$ which satisfies the energy inequality (\ref{ei}).}
\end{Theorem}

\bRemark{R1}
\begin{enumerate}
\item Theorem \ref{TM1} still holds provided one considers in the momentum equation at its right hand side the term $\vr\vc f$, $\vc f\in L^\infty(Q_T)$, corresponding
to the action of large external forces. The necessary changes in the weak formulation and in the relative energy inequality  in order to accommodate the presence of this term are left to the reader.
\item Conditions on the regularity $p$, $\vr_B$ and $\vu_B$ in Theorem \ref{TM1} could be slightly weakened, up to
${p}$ continuous on $[0,\infty)$, locally Lipschitz on $[0,\infty)$, $\vr_B\in L^\infty(\partial\Omega)$, $\vu_B\in W^{1,\infty}(\partial\Omega)$, at expense of some additional technical difficulties.
\item We shall perform the proof in all details in the case $d=3$ assuming tacitly that both $\Gamma_{\rm in}$ and $\Gamma_{\rm out}$ have non zero
$(d-1)$-Hausdorff measure. Other cases, namely the case  $d=2$ is left to the reader as an exercise.
\end{enumerate}
\eR


Our second main result is the following.

\bTheorem{WS1}{\rm [Stability and Weak-strong uniqueness principle]}
Let $\Omega$ be a bounded Lipschitz domain.
 Suppose that the pressure is in addition to (\ref{pressure}), twice continuously differentiable on $(0,\infty)$ and obeys
 \begin{equation}\label{w1+}
c(\vr+ H(\vr))\ge p(\vr)\; \mbox{for all $\vr\ge\overline R$},
\end{equation}
where $\overline R$, $c$ are some positive constants.
Assume that  the initial data
$\vr_0$, $\vc u_0$ verify condition (\ref{feid}) and boundary data  satisfy condition (\ref{M1}). Let $\vu_\infty$ be an Lipschitz extension of $\vu_B$ satisfying conditions (\ref{Le}).

Let $(\vr,\vc u)$ be a dissipative weak solution to the Navier-Stokes equations (\ref{i1}-\ref{i5}) associated to the extension $\vu_\infty$
of $\vu_B$. Let $(r,\vc U)$
that belongs to the class
\begin{equation}\label{s1}
0<\underline r\le r\le\overline r<\infty;\quad\vc U \in L^\infty((0,T)\times\Omega)
\end{equation}
$$
\partial_t r,\partial_t\vc U,\Grad r,\Grad\vc U\in L^2(0,T;C(\overline\Omega))\footnote{The requested regularity of some of this derivatives
can be slightly weaken by a more detailed (elementary) analysis as in \cite{FeJiNo}.}
$$
be a strong solution of the same equations with  initial data $(r_0,\vc U_0)$ and boundary data $(r|_{\Gamma_{\rm in}}=r_B,\vc U|_{\partial\Omega}=
\vu_B).$

\begin{enumerate}
\item
Then there exists
$$
c= c\Big(\mu, T, |\Omega|, {\rm diam}\Omega, \mu, \underline{\mathfrak{r}}, \overline{\mathfrak{r}},
|p,p',H'|_{C^1([\underline{\mathfrak{ r}}/2, 2\overline{\mathfrak{ r}}])},
\|\partial_t\vc V, \nabla\vc V, \nabla r\|_{L^2(0,T;L^\infty(\Omega;R^{15}))}\Big)>0
$$
(where  $\underline{\mathfrak r}=\min\{\underline r,\underline\vr_B\}$,
$\overline{\mathfrak r}=\max\{\overline r,\overline\vr_B\}$),
such that
\bFormula{stability1+}
{\cal E}(\vr,\vc v| r ,\vc V)(\tau)\le c\Big( {\cal E}(\vr_0,\vc v_0| r_0 ,\vc U_0) +
\|\vr_B- r_B\|_{L^1(\Gamma_{\rm in})}\Big)
\eF
for a.e. $\tau\in (0,T)$, where $\vc v=\vu-\vu_\infty$, $\vc V=\vc U-\vu_\infty$, $\vc v_0=\vc u_0-\vc u_\infty$,
$\vc V_0=\vc U_0-\vc u_\infty$. In the above, we have denoted
\bFormula{refunct}
{\cal E}(\vr,\vc v| r ,\vc V):=
\intO{\Big(\frac 12\vr|\vc v-\vc V|^2+E(\vr|r)\Big)}
\eF
the relative energy functional.
\item
In particular, if $(\vr_0,\vu_0)=(r_0,\vc U_0)$, and $r_B=\vr_B$ then
$$
\vr=r,\; \vc u=\vc U\; \mbox{in $[0,T]\times\Omega$}.
$$
\end{enumerate}
\eT
\noindent

\bRemark{RWS}
\begin{enumerate}
\item  The (dissipative) weak-strong uniqueness holds for pressure functions fairly beyond the conditions (\ref{pressure1})
guaranteeing existence of weak solutions. An example of an admissible class of pressure functions is the class
$p\in C[0,\infty)\cap C^2(0,\infty)$, $p(0)=0$, $p'(\vr)>0$ and
\begin{equation}\label{dodano}
0<\frac 1{p_\infty}\le\liminf_{\vr\to\infty}\frac {p(\vr)}{\vr^\gamma}\le\limsup_{\vr\to\infty}\frac {p(\vr)}{\vr^\gamma}\le p_\infty<\infty,\;\gamma\ge 1.
\end{equation}
Indeed, one can easy verify, that  condition (\ref{dodano}) yields condition (\ref{w1+}) from Theorem \ref{TWS1}. In particular, condition (\ref{pressure1})
implies (\ref{dodano}), whence also (\ref{w1+}). Consequently, dissipative weak solutions constructed in Theorem \ref{TM1} satisfy the weak-strong
uniqueness principle.
\item  Existence of strong solutions at least on a short time interval is well known.
Here we report the following existence result of Valli and Zajaczkowski \cite[Theorem 2.5]{VAZA}.
\bLemma{S1}
Let $D$ be a positive constant, $\Omega$  a bounded domain  of class $C^3$ and $p\in C^2((0,\infty))$.
Let
$$
\vu_\infty\in W^{3,2}(\Omega),\; \vu_B=\vu_\infty|_{\partial\Omega},\; \vr_B\in W^{2,2}(\Gamma_{\rm in})
$$
where the inflow boundary $\Gamma_{\rm in}$ is define in (\ref{i5}).
 Assume further that
\begin{equation}\label{itar1f}
\vc u_{0}\in W^{3,2}(\Omega;R^3),\;\vu_0|_{\partial\Omega}=\vu_B
\end{equation}
$$
\vr_0\in W^{2,2}(\Omega),\;\inf_{\Omega} \vr_0>0, \; \vr_0|_{\Gamma_{\rm in}}=\vr_B,
$$
where
\begin{equation}\label{itar1+f}
{\rm div}(\vr_0\vu_0)|_{\Gamma_{\rm in}}=0,
\end{equation}
$$
\frac 1{\vr_0}\Big(-\Grad p(\vr_0)+{\rm div}\tn S(\Grad\vc u_0)-\vr_0\vc u_0\cdot\Grad\vc u_0\Big)\Big|_{\partial\Omega}=0.
$$

Then there exists $T=T_0(D)$ such that if
\begin{equation}\label{itar2}
\|\vr_0\|_{ W^{2,2}(\Omega)}+\|\vc u_0\|_{ W^{3,2}(\Omega)}+1/\inf_{\Omega} \vr_0 +\|\vu_\infty\|_{W^{3,2}(\Omega)}+
\|\vr_B\|_{W^{2,2}(\Gamma_{\rm in})}\le D,
\end{equation}
then the problem (\ref{i1}--\ref{i5}) admits a unique strong solution (in the sense a.e. in $(0,T)\times\Omega$) in the class
\begin{equation}\label{itar3}
\vr\in C([0,T);W^{2,2}(\Omega)),\; \vc u-\vu_\infty\in C([0,T);W^{3,2}\cap W_0^{1,2}(\Omega;R^3))\cap L^2(0,T; W^{4,2}(\Omega;R^3)),
\end{equation}
$$
\partial_t \vr\in C([0,T);W^{1,2}(\Omega)),\;\partial_t\vc u\in L^2([0,T);W^{2,2}(\Omega;R^3)).
$$
In particular,
\begin{equation}\label{itar4}
0<\underline r\equiv\inf_{(t,x)\in (0,T)\times\Omega} \vr(t,x)\le \sup_{(t,x)\in (0,T)\times\Omega} \vr (t,x)\equiv\overline r<\infty.
\end{equation}
\eL
\end{enumerate}
\eR

\section{Approximate problem}
\label{A}
Our goal is to construct a dissipative weak solution, the existence of which is claimed in Theorem \ref{TM1}. We shall use the approximating procedure
suggests in \cite{JiNoSIMA}. We add to the pressure an artificial pressure (small parameter $\delta>0$), regularize the continuity equation by adding an
artificial viscosity term to equation (\ref{i1}) endowing at the same time the momentum equation by a compensation term in order to keep
the energy inequality - this is so far standard - and add to the momentum equation a monotone dissipation  (the latter
adjustment is mostly technical) - all these three ingredients parametrized by a small parameter $\ep>0$. Moreover, the regularized continuity equation is completed with boundary conditions through a boundary operator which gives in the limit
$\ep\to 0$ the general inflow/outflow conditions.

The  approximating system of equations reads:
\begin{equation}
\label{A1}
\partial_t\vr-\ep \Del \vr  + \Div (\vr \vu) = 0,
\end{equation}
\begin{equation}
\label{A2}
\vr(0,x)=\vr_0(x),\;\left( -\ep \Grad \vr\cdot\vc n + \vr v \right)|_{I\times\partial \Omega} = g,
\end{equation}
where
$$
\left\{\begin{array}{c}
\vu_B\cdot\vc n\;\mbox{on $\Gamma_{\rm in}$},\\
0\;\mbox{on $\partial\Omega\setminus\Gamma_{\rm in}$}
\end{array}\right\}\equiv v,\quad \vr_B v\equiv g.
$$
\begin{equation}
\label{A3}
\partial_t(\vr\vu)+
\Div (\vr \vu \otimes \vu) + \Grad p_{\delta} (\vr) = \Div \mathbb{S}(\Grad \vu) -\ep \Grad\vr\cdot\Grad\vu   +{ \ep{\rm div}\Big(|\Grad(\vu-\vu_\infty)|^2\Grad (\vu-\vu_\infty)\Big)}
\end{equation}
\begin{equation}
\label{A4}
\vu(0,x)=\vu_0(x),\;
\vu|_{I\times\partial \Omega} = \vu_B,
\end{equation}
with positive parameters $\ep > 0$, $\delta > 0$,
where
we have denoted
\bFormula{pdelta}
p_\delta(\vr) = p(\vr)+\delta\vr^\beta,\quad \beta> \max\{\gamma,9/2\}
\eF
and where $\vu_\infty$ is an extension of $\vu_B$ defined in (\ref{Le}).\footnote{The exact choice of $\beta$ is irrelevant from the point of view of the final result provided it is sufficiently large.}
Next, following \cite{JiNoSIMA}, we shall define the generalized solutions to the approximate problem (\ref{A1}--\ref{A4}).\footnote{The only difference with respect
to \cite{JiNoSIMA} in this definition is the fact, that the test function in the equation (\ref{Aw1}) does not vanish  at the outflow boundary and the energy
inequality is more precise containing also all boundary terms. This is essential for the construction of dissipative solutions.}
\\ \\
{\bf Definition \ref{A}.1}
{\it Let $\vu_\infty\in W^{1,\infty}(\Omega;\R^3)$ be a Lipschitz extension of $\vu_B$ staisfying
(\ref{hatU})\footnote{For bounded Lipschitz domains such extension always eists, cf. Remark \ref{RR1-}}.
A couple $(\vre,\vue)$ and associated tensor field $\tn Z_\ep$ is a generalized solution of the sequence of problems (\ref{A1}--\ref{A4})$_{\ep>0}$
iff the following holds:
\begin{description}
\item{1.} It belongs to the functional spaces:
\bFormula{fsw}
\vre \in L^\infty (0,T; L^\beta(\Omega))\cap L^2(0,T; W^{1,2}(\Omega)), \ 0 \leq \vre \ \mbox{a.e. in}\ (0,T)\times \Omega,
\eF
$$
\vc v_\ep:=\vu_\ep-\vu_\infty \in L^2(0,T; W_0^{1,2}(\Omega; R^3))\cap L^4(0,T;W_0^{1,4}(\Omega; R^3)),
$$
$$
\tn Z_\ep\to 0\;\mbox{in  $L^{4/3}(Q_T;\R^3)$ as $\ep\to 0$}.
$$

\item{2.}
Function $\vre\in C_{\rm weak}([0,T], L^\beta(\Omega))$ and the integral identity
\begin{align}\label{Aw1}
&\intO{\vre(\tau,x)\varphi(\tau,x)} - \intO{\vr_0(x)\varphi(0,x)}\nonumber\\
&\qquad=
\int_0^\tau\intO{\Big(\vre\partial_t\varphi + \vre \vue \cdot \Grad \varphi-\ep\Grad\vre\cdot\Grad\varphi\Big) }{\rm d}t
\nonumber\\
&\qquad\qquad
- \int_0^\tau\int_{\Gamma_{\rm in}} \vr_B \vu_B \cdot \vc{n} \varphi \ {\rm d}S_x{\rm d} t
- \int_0^\tau\int_{\Gamma_{\rm out}} \vre \vu_B \cdot \vc{n} \varphi \ {\rm d}S_x{\rm d} t
\end{align}
holds for any $\tau\in [0,T]$ and $\varphi \in C_c^1([0,T]\times\overline\Omega)$.
\item{3.} Function
$\vre\vue\in  C_{\rm weak}([0,T], L^{\frac{2\beta}{\beta +1}}(\Omega;\R^3))$,  and the integral identity
\begin{align}\label{Aw2}
&\intO{\vre\vc v_\ep(\tau,\cdot)\cdot\bfphi(\tau,\cdot)} - \intO{\vr_0\vc v_0(\cdot)\bfphi(0,\cdot)}= { -\int_0^\tau\intO{\tn Z_\ep:\Grad\varphi}{\rm d} t}\nonumber\\
&
\qquad+\int_0^\tau\int_\Omega  \Big(\vre\vc v_\ep\partial_t\varphi +\vre \vue \cdot \Grad \bfphi\cdot\vue
-\vre\vue\cdot\Grad(\vu_\infty\cdot \bfphi)
+ p_\delta(\vre) \Div \bfphi\nonumber\\
&
\qquad\qquad-\mathbb{S}(\Grad \vue) : \Grad \bfphi
+\ep\Grad\vre\cdot\Grad(\vu_\infty\cdot\bfphi)
-
\ep\Grad\vre\cdot\Grad\vue\cdot\bfphi\Big) {\rm d} x{\rm d}t
\end{align}
holds for any  $\tau\in [0,T]$ and any $\bfphi \in C^1_c ([0,T]\times\Omega; R^3)$.
\item{4.} Energy inequality
\begin{align}\label{Aw3}
&\intO{\Big(\frac 12\vre|\vc v_\ep|^2
+H_\delta(\vre)\Big)}
+
\int_0^\tau\int_{ \Gamma_{\rm in}} E_\delta(\vr_B|\vre)|\vu_B\cdot\vc n|{\rm d}S_x{\rm d}t+
\int_0^\tau\int_{ \Gamma_{\rm out}} H(\vre)|\vu_B\cdot\vc n|{\rm d}S_x{\rm d}t\nonumber\\
&\qquad
+\int_0^\tau\intO{\Big(\tn S(\Grad\vue):\Grad\vc v_\ep
+\ep H_\delta''(\vre)|\Grad\vre|^2 { +\ep|\Grad\vc v_\ep|^4}\Big)}{\rm d}t\nonumber\\
&
\le
{ \intO{\Big(\frac 12\vr_0|\vu_0-\vu_\infty|^2+H_\delta(\vr_0)\Big)}}
- \int_0^\tau\int_{\Gamma_{\rm in}} H_\delta(\vr_B)\vu_B\cdot\vc n{\rm d}S_x{\rm d} t\nonumber\\
&
\qquad-\int_0^\tau\intO{\Big(p_\delta(\vre){\rm div}\vu_\infty
+\vre\vue\cdot\Grad\vu_\infty\cdot\vc v_\ep-\ep\Grad\vr\cdot\Grad\vc v_\ep\cdot\vu_\infty\Big)}{\rm d}t
\end{align}
holds for a.e. $\tau\in (0,T)$.  In the above,
\bFormula{Hfd}
H_\delta(\vr)= H(\vr) + \delta H^{(\beta)}(\vr),\; H^{(\beta)}(\vr)=
\vr\int_1^\vr z^{\beta-2} {\rm d}z=\frac 1{\beta-1}\vr^\beta
\eF
and
\begin{align}\label{Ed}
& E_\delta(\vr|r)= {E}(\vr|r)+\delta{E}^{(\beta)}(\vr|r),~
E^{(\beta)}(\vr|r)= H^{(\beta)}(\vr)-[H^{(\beta)}]'(r)(\vr-r)-H^{(\beta)}(r).
\end{align}
\end{description}
}

The starting point in the construction of the dissipative solutions for system (\ref{i1}--\ref{i5}) is the existence theorem for the generalized approximating problem (\ref{A1}--\ref{A4}). It is announced in the next lemma. \footnote{The only difference with respect
to \cite[Lemma 4.1]{JiNoSIMA} in Lemma \ref{LLemmaA} is the fact, that the test function in equation (\ref{Aw1}) do not vanish  at the outflow boundary and the energy
inequality is more precise containing also all boundary terms. This is essential for the construction of dissipative solutions.}

\bLemma{LemmaA}
Let $\Omega $ be a domain of class { $C^{2}$}. Let $(\vr_B,\vu_B)$ verify assumptions  (\ref{M1}) and let initial  and boundary data
verify
\bFormula{idod}
\vu_0\in L^2(\Omega),\quad \vr_0\in W^{1,2}(\Omega),\;0<\underline\vr\le\vr_0\le\overline\vr<\infty,
\eF
\bFormula{bdod}
0<\underline\vr_B\le\vr_B\le\overline\vr_B<\infty.
\eF
Then for any continuous extension $\vu_\infty$ of $\vu_B$ in class (\ref{Le})  there exists a  generalized solution
$(\vr_\epsilon,\vu_\epsilon )$ and $\tn Z_\ep$ to the sequence of approximate problems
(\ref{A1} - \ref{A4})$_{\ep\in (0,1)}$ - which belongs to the functional spaces (\ref{fsw}), satisfies the weak formulations (\ref{Aw1}--\ref{Aw2})
and verifies the energy inequality
(\ref{Aw3}) -  with the following extra properties:
\begin{description}
\item {\rm (i)} In addition to (\ref{fsw}) it belongs to functional spaces:
\begin{equation}\label{ts5.23}
 \vr_\ep\in L^{\frac 53 \beta}(Q_T), \;\sqrt\vre, \vr_\epsilon^{\frac \beta 2}\in L^2(I,
W^{1,2}(\Omega )),\;
\partial_t\vr_\epsilon\in L^{4/3}(Q_T),
\nabla^2\vr_\epsilon \in L^{4/3}
(Q_T).
\end{equation}
\item{\rm (ii)} In addition to the weak formulation (\ref{Aw1}), the couple $(\vr_\ep,\vu_\ep)$ satisfies the equation (\ref{Aw1}) in the strong sense,
meaning, it verifies equation (\ref{A1}) with $(\vr_\ep,\vu_\ep)$
a.e. in $Q_T$, boundary identity (\ref{A2}) with $(\vr_\ep,\vu_\ep)$ a.e. in $(0,T)\times \partial\Omega$ and initial conditions in the sense
$\lim_{t\to 0+}\|\vre(t)-\vr_0\|_{L^{4/3}(\Omega)}=0$.
\item{\rm (iii) } The couple $(\vr_\ep,\vu_\ep)$ satisfies identity
\begin{equation} \label{A11}
\partial_t b(\vr_\ep)+\ep b''(\vr_\ep) |\Grad \vr_\ep|^2  - \ep \Div (b'(\vr_\ep) \Grad \vr_\ep) + \Div (b(\vr_\ep) \vu_\ep)+
\left[ b'(\vr_\ep) \vr_\ep - { b(\vr_\ep)} \right] \Div \vu_\ep = 0
\end{equation}
a.e. in $(0,T)\times\Omega$  with any $b\in C^2[0,\infty)$, where the space-time derivatives have to be understood in the sense a.e.
\end{description}
\eL
\begin{Remark}\label{R1+e}
Identity (\ref{A11}) holds in the weak sense
$$
\intO{b(\vre(\tau))\varphi(\tau)}-\intO{b(\vr_0)\varphi(0)}=
\int_0^\tau\int_{\partial\Omega}\Big(\ep b'(\vre)\Grad\vre- b(\vre)\vue\Big)\cdot\vc n{\rm d}S_x{\rm d}t
$$
$$
+\int_0^\tau\intO{\Big[b(\vre)\partial_t\varphi+(b(\vre)\vue-\ep b'(\vre)\Grad\vre)\cdot
\Grad\varphi -\varphi\Big(\ep b''(\vre)|\Grad\vre|^2+(\vre b'(\vre)-b(\vre)){\rm div}\vue\Big)\Big]}{\rm d}t
$$
with any $\tau\in [0,T]$ and $\varphi\in C^1_c([0,T]\times\overline\Omega)$ with any
$b$ whose growth (and that one of its derivatives) in combination with (\ref{ts5.23}) guarantees
$b(\vre)\in C_{\rm weak}([0,T]; L^1(\Omega))$, existence of traces and integrability of all terms appearing at the r.h.s.
\end{Remark}

\section{Construction of the generalized solutions to the approximate problem}\label{4}

We recall here the main building blocks of the construction of generalized solutions to the approximate problem (\ref{A1}--\ref{A4}).
We do not intend to describe the whole process in all details, since it is available in \cite[Section 4]{JiNoSIMA}, but only
its main parts.

\subsection{Galerkin type approximation and energy inequality}

1. The first building block in the construction of generalized solutions to the approximate problem (\ref{A1}--\ref{A4}) is the following theorem dealing
with
the parabolic problem
(\ref{A1}--\ref{A2}). It reads, cf. \cite[Lemma 4.3]{JiNoSIMA}:
\bLemma{Lparabolic}
Suppose that $\Omega$ is a bounded domain of class $C^{2}$ and assume further that
{ $\vr_0\in W^{1,2}(\Omega)$},  { $\vu|_{(0,T)\times\partial\Omega}=\vu_B$}, $v, g\in C^1(\partial\Omega)$. Then we have:
\begin{enumerate}
\item
The parabolic  problem (\ref{A1}--\ref{A2}) admits for any $\vc u \in L^\infty(0,T; W^{1,\infty}(\Omega))$, a unique solution
$\vr = S(\vc u)$
in the class
\bFormula{parabolicr}
\vr\in L^2(0,T;W^{2,2}(\Omega))\cap W^{1,2}(0,T;L^2(\Omega)).
\eF
\item Suppose that
$$
\underline\vr\le\vr_0(x) \le\overline\vr\;\mbox{for a. a. $x\in\Omega$},\quad \underline\vr\le\vr_B(x) \le\overline\vr\;\mbox{for all $x\in\Gamma_{\rm in}$}.
$$
Then
\bFormula{par4}
\underline\vr {\rm exp}\Big(-\int_0^\tau\|{\rm div}\vc u(s)\|_{L^\infty(\Omega)}{\rm d}s\Big)\le\vr(\tau,x)\le  \overline\vr
{\rm exp}\Big(\int_0^\tau\|{\rm div}\vc u(s)\|_{L^\infty(\Omega)}{\rm d}s\Big),
\eF
in particular,
$$
\underline\vr e^{-K\tau}\le\vr(t,x)\le \overline\vr e^{K\tau}
$$
for all $\tau\in[0,T]$ provided $\vc u$ verifies condition
\bFormula{K}
\|\vu\|_{L^\infty(Q_T)}+ \|{\rm div}\vu\|_{L^\infty(Q_T)}\le K.
\eF
\end{enumerate}
\eL
\vspace{5mm}
\noindent
2. The second building block in the construction of generalized solutions to the approximate problem (\ref{A1}--\ref{A4}) is the Galekin approximation
of the momentum equation (\ref{A3}--\ref{A4}).
In order to write it down, we introduce 
\bFormula{X}
X={\rm span}\{\Phi_i\}_{i=1}^N\;\mbox{where ${\cal B}:=\{\Phi_i\in C^\infty_c(\Omega)\,|\, i\in \tn N^*\}$ is an orthonormal basis in $L^2(\Omega;\R^3)\}$},
\eF
a finite dimensional real Hilbert space with scalar product $(\cdot,\cdot)_X$ induced by the scalar product in $L^2(\Omega;R^3)$ and
$\|\cdot\|_X$ the norm induced by this scalar product. We denote by $P_N$ the orthogonal projection of $L^2(\Omega;R^3)$ to $X$.

The Galerkin approximation of system (\ref{A1}-\ref{A4}) reads: {\it Given $\vr_B$, $\vu_\infty$  (whence also $\vu_B$), and $\vu_0=\vc v_0-\vu_\infty$, $\vr_0$,
find the couple $(\vr,\vu)$, $\vu=\vc v+\vu_\infty$,
\begin{equation}\label{clv}
\vc v \in C(\overline I; X), \, \partial_t\vc v\in L^2(I; X),\; 0\le\vr\in L^2(I;W^{2,2}(\Omega))\cap W^{1,2}( I;L^2(\Omega))
\end{equation}
such that
\begin{align}\label{gal}
&\intO{\vr\vc v(t)\cdot\Phi}-\intO{\vr_0\vc v_0\Phi}=\int_0^t\intO{\Big({\rm div}{\tn S}(\Grad\vu){
+\ep{\rm div}\Big(|\Grad(\vu-\vu_\infty)|^2\Grad (\vu-\vu_\infty)\Big)}
\nonumber\\
&
\qquad-\Grad p_\delta(\vr) -{\rm div}(\vr\vu\otimes\vu)
-\ep\Grad\vr\cdot\Grad\vu{-\partial_t\vr\vu_\infty }\Big)\cdot\Phi}{\rm d}t
\end{align}
holds for all $t\in [0,T]$ with any $\Phi\in X$, where $\vc u=\vc v+\vu_\infty$ and
\begin{equation}\label{ap+}
\mbox{$\vr$ solves parabolic problem (\ref{A1}--\ref{A2}).}
\end{equation}
}

The next lemma ensures existence of solutions to problem (\ref{clv}--\ref{ap+}).
\bLemma{LGalerkin}
Let $\Omega$, $(\vr_0,\vu_0)$, $(\vr_B,\vu_B)$ verify assumptions of Lemma \ref{LLemmaA}.
\begin{enumerate}
\item
Then for any continuous extension $\vu_\infty$ of $\vu_B$ in class (\ref{Le})  the Galerkin type problem
(\ref{gal}--\ref{ap+}) admits a unique solution $(\vr,\vc v)$ in class (\ref{clv}). Moreover,
there are constants $0<\underline c=\underline c(N)<\overline c=\overline c(N)$ such that
\begin{equation}\label{rho>0}
\forall t\in \overline I, \; \underline c\le\vr(t,\cdot)\le\overline c\;\mbox{a.e. in $\Omega$}.
\end{equation}
\item The solution satisfies energy inequality:
\begin{align}\label{galen}
&\intO{\Big(\frac 12\vr|\vc v|^2+H_\delta(\vr)\Big)(\tau)} +\int_0^\tau\int_{\Gamma_{\rm in}} E_\delta(\vr_B|\vr)|\vu_B\cdot\vc n|{\rm d}S_x{\rm d}t
+\int_0^\tau\int_{\Gamma_{\rm out}} H_\delta(\vr)|\vu_B\cdot\vc n|{\rm d}S_x{\rm d}t
\nonumber\\
& \qquad+\ep\int_0^\tau\intO{|\Grad \vc v|^4}{\rm d}t
+ \ep\int_0^\tau\intO{H''_\delta(\vr)|\Grad\vr|^2}{\rm d}t
+\int_0^\tau\intO{\tn S({\Grad\vc u}):\Grad\vc v}{\rm d}t
\nonumber\\
&\le  \intO{\Big(\frac 12\vr_0|\vc v_0|^2+H_\delta(\vr_0)\Big)}
-\int_0^\tau\int_{\Gamma_{\rm in}} H_\delta(\vr_B)\vu_B\cdot\vc n{\rm d}S_x{\rm d}t
\nonumber\\
&\qquad+
\int_0^\tau\intO{\Big(-p_\delta(\vr){\rm div}\vu_\infty
-\vr\vu\cdot\Grad\vu_\infty\cdot\vc v
+\ep\Grad\vr\cdot\Grad\vc v\cdot\vu_\infty\Big)}{\rm d}t,
\end{align}
where $H_\delta$ and $E_\delta$ are defined in (\ref{Hfd}) and (\ref{Ed}).
\end{enumerate}
\eL
\noindent
{\bf Proof of Lemma \ref{LLGalerkin}}\\
We refer to \cite[Section 4.3]{JiNoSIMA} for the proof of the existence part (Item 1.). Energy inequality (\ref{galen}) - Item 2. -
is derived in formula (4.40) of the same paper. We repeat its proof here for the sake of completeness, since
the energy inequality  plays crucial role in the genesis of the
proof of the exact energy inequality (\ref{Aw3}) and consequently in the proof of the existence of
dissipative solutions.

We have at our disposal the couple $(\vr=\vr_N=S(\vu_N),\vu=\vu_N)$ a solution of problem (\ref{gal}--\ref{ap+}) in the class (\ref{clv}) and (\ref{rho>0}).
We are going to derive for this couple the energy inequality (\ref{galen}).

First, we multiply equation (\ref{ap+}) by $H_\delta'(\vr)$ and integrate over $\Omega$  to deduce
\begin{align}\label{rena}
&\partial_t\intO{ H_\delta(\vr)}+\ep\intO{H_\delta''(\vr)|\Grad\vr|^2}\nonumber\\
&\qquad+\int_{\Gamma_{\rm in}}\Big(H_\delta(\vr_B) -H'_\delta(\vr){(\vr_B-\vr)}
-H_\delta(\vr)\Big)|v|{\rm d}S_x +\int_{\Gamma_{\rm out}}H_\delta(\vr)\vu_B\cdot\vc n{\rm d}S_x
\nonumber\\
&=-\intO{p_\delta(\vr){\rm div}\vu} - \int_{\Gamma_{\rm in}} H_\delta(\vr_B)\vu_B\cdot\vc n {\rm d}S_x,
\end{align}
where we have used several times integration by parts, take into account boundary conditions (\ref{ap+}), and where
$\vc v=\vu-\vu_\infty$.

Further, we deduce from (\ref{gal})
$$
\int_0^\tau\intO{\Big(\partial_t(\vr \vc v)\cdot\vc v -\vr\vu\otimes\vu:\Grad\vc v\Big)}{\rm d}t + 
\int_0^\tau\intO{\Big(\tn S(\Grad\vu):\Grad\vc v+\ep|\Grad \vc v|^4\Big)}{\rm d}t
$$
$$
-\int_0^\tau\intO{p_\delta(\vr){\rm div}\vc v}{\rm d}t + \int_0^\tau\intO{\ep \Grad\vr\cdot\Grad\vu\cdot\vc v}{\rm d}t
+\int_0^\tau\intO{\partial_t\vr\vc u_\infty\cdot\vc v}{\rm d}t=0,
$$
where
by virtue of (\ref{ap+}) (after several integrations by parts and recalling that $\vu=\vu_\infty+\vc v$),
\begin{align*}
&\intO{\Big(\partial_t(\vr\vc v)\cdot\vc v -\vr\vu\otimes\vu:\Grad \vc v\Big)} \\
&\qquad\qquad=\intO{\Big(\partial_t\vr\vc v^2+\frac 12\vr\partial_t(\vc v^2)+
\frac 12{\rm div}(\vr\vu)\vc v^2 -\vr\vu\cdot\Grad\vc v\cdot\vu_\infty\Big)}
\\
&\qquad\qquad
=\intO{\Big(\frac 12\partial_t(\vr\vc v^2)
-\ep \Grad\vr\cdot\Grad \vc v\cdot\vc v
-\vr\vu\cdot\Grad\vc v\cdot\vu_\infty\Big)}
\end{align*}
and 
$$
\intO{\partial_t\vr\vu_\infty\cdot\vc v}=
\intO{\Big(\vr\vu\cdot\Grad\vc v\cdot\vu_\infty+\vr\vu\cdot\Grad\vu_\infty\cdot\vc v-
\ep\Grad\vr\cdot\Grad\vu_\infty\cdot \vc v-\ep\Grad\vr\cdot\Grad \vc v\cdot\vu_\infty\Big)}.
$$
Therefore,
$$
\int_0^\tau\intO{\partial_t(\vr \vc v)\cdot\vc v -\vr\vu\otimes\vu:\Grad\vc v+
\ep \Grad\vr\cdot\Grad\vu\cdot\vc v + \partial_t\vr\vu_\infty\cdot\vc v\Big)}
$$
$$
=\int_0^\tau\intO{\Big(\vr\vu\cdot\Grad\vu_\infty\cdot\vc v -\ep\Grad\vr\cdot\Grad\vc v\cdot\vu_\infty\Big)}{\rm d}t.
$$
This together with (\ref{rena}) yields inequality (\ref{galen}). Lemma \ref{LLGalerkin} is thus proved.

\subsection{Uniform bounds with respect to $N$ and limit $N\to\infty$}

In order to get uniform bounds for the sequence $(\vr_N,\vu_N)$ we still need the conservation of mass
\bFormula{mi+}
\intO{\vr(\tau)}+\int_0^\tau\int_{\Gamma_{\rm out}}\vr\vu_B\cdot\vc n{\rm d}S_x{\rm d}t
=\intO{\vr_0} + \int_0^\tau\int_{\Gamma_{\rm in}}\vr|\vu_B\cdot\vc n|{\rm d}S_x{\rm d}t
\eF
(which we obtain from (\ref{ap+})). We also need to 
to test
equation (\ref{A1})$_{(\vr_N,\vu_N)}$ by $\vr_N$ in order to get, after
several integrations by parts,
\begin{align}\label{Gr1}
&{
\frac 12\intO{\vr^2(\tau)}+\frac 12\int_0^\tau \int_{\partial\Omega}\vr^2 |\vu_B\cdot\vc n|{\rm d}S_x{\rm d}t +\ep\int_0^\tau\intO{|\Grad\vr|^2}{\rm d}t
}
\nonumber\\
&\qquad
=\frac 12\intO{\vr_0^2} + \int_0^\tau \int_{\Gamma_{\rm in}}\vr\vr_B |\vu_B\cdot\vc n|{\rm d}S_x{\rm d}t-\frac 12\int_0^\tau\int_\Omega\vr^2{\rm div}\vu{\rm d}x
{\rm d}t.
\end{align}

{Recalling structural assumptions (\ref{pressure}), (\ref{pressure1})  for $p$, definitions  (\ref{H}), (\ref{pdelta}),(\ref{Hfd}),
we deduce from the energy inequality (\ref{galen}) and inequalities (\ref{mi+}--\ref{Gr1})
 the following uniform bounds with respect to $N$ for the sequence $(\vr_N=S(\vu_N),\vu_N=\vu_\infty+\vc v_N)$ of Galerkin
solutions to the problem (\ref{gal}--\ref{ap+}):
\begin{align}
&\|\vr_N|\vu_N|^2\|_{L^\infty(I,L^1(\Omega ))}\le
 L({\rm data}),\label{ts1?}\\
&
\|\vu_N\|_{L^2(I,W^{1,2}(\Omega ))}\le L({\rm data}),\label{ts2?}\\
&
\|\vr_N\|_{L^\infty(I,L^\gamma(\Omega))}\le
L({\rm data}),\label{ts2?+}\\
&
\delta^{1/\beta} \|\vr_N\|_{L^\infty(I,L^\beta(\Omega))}\le
L({\rm data}),\label{ts3?}\\
&
\ep\|\nabla{\vr_N}\|^2_{L^2(Q_T)} + \ep\|\nabla(\vr_N^{\beta/2})\|^2_{L^2(Q_T)} \le
L({\rm data},\delta),\label{ts4?}\\
&
\ep^{1/\beta}\|\vr \|_{L^\beta((0,T)\times\partial\Omega)}\le  L({\rm data},\delta),
\label{nova+}\\
&
\ep\|\vu_N-\vu_\infty\|^4_{L^4(0,T; W^{1,4}(\Omega))}\le  L({\rm data},\delta).\label{vu4}
\end{align}
In the above and hereafter
$$
\mbox{"data" stands for}\; \intO{\Big(\frac 12\vr_0|\vc v_0|^2+H(\vr_0)\Big)},\,\|\vu_\infty\|_{W^{1,\infty}(\Omega)},\,\underline \vr,\,\overline\vr,\,
\underline\vr_B,\,\overline\vr_B,\,\underline H=\inf_{\vr>0}H(\vr).
$$
Due to the above estimates, expression  ${\rm div}(\vr_N\vu_N)$ is bounded in $L^{4/3}(Q_T)$. We can thus return to equation (\ref{A1}--\ref{A2})
-with $(\vr_N,\vu_N)$- and consider it as parabolic problem with operator $\partial_t\vr-\ep\Delta\vr$ in $(0,T)\times\Omega$ with right
hand side $-{\rm div} (\vr_N\vu_N)$, and boundary operator $-\ep\vc n\cdot\Grad\vr +v \vr$ in $(0,T)\times\partial\Omega$ with right hand side $\vr_B v$.
The maximal parabolic regularity theory, see e.g. \cite[Theorem 2.1]{DHP}, yields
that
\begin{equation}\label{ts8?}
\|\partial_t\vr_N\|_{L^{4/3}(Q_T)} + \|\vr_N\|_{L^{4/3}(0,T;W^{2,4/3}(\Omega))}\le L({\rm data},\delta,\ep ).
\end{equation}

Estimates (\ref{ts1?}--\ref{ts8?}) yield, in particular, existence of a subsequence $(\vr_N,\vu_N)$ not relabeled such that
\begin{align}
&\vr_N\rightharpoonup\vr\;\mbox{ in $L^{4/3}(0,T;W^{2,4/3}(\Omega))$ and in $L^2(0,T;W^{1,2}(\Omega))$},\quad \partial_t\vr_N\rightharpoonup\partial_t\vr\;
\mbox{in $L^{4/3}(Q_T)$}\label{zac}\\
&
\vr_N\rightharpoonup\vr\;\mbox{in $L^\beta((0,T)\times\partial\Omega)$},\label{nova+++}\\
&
\vu_N\rightharpoonup\vu\;\mbox{ (weakly) in $L^4(0,T;W^{1,4}(\Omega))$},\label{cN3}\\
&
\ep|\Grad(\vu_N-\vu_\infty)|^2\Grad(\vu_N-\vu_\infty)\rightharpoonup \tn Z\equiv\tn Z_\ep \;\mbox{weakly in $L^{4/3}(Q_T; R^9)$,}\label{cN4}
\end{align}
 where
$$
\|\tn Z_\ep\|_{L^{4/3}(Q_T)}\to 0\;\mbox{as $\ep\to 0$}.
$$

This information is enough to show that $(\vr,\vu, \tn{Z})$ to class (\ref{fsw}),(\ref{ts5.23}) and allows to pass to the limit $N\to\infty$ in the system (\ref{gal}--\ref{ap+}) and in the inequality (\ref{galen}) by using only the standard compactness arguments (which are Sobolev embeddings, Arzela-Ascoli theorem and Lions-Aubin Lemma). In particular, relation (\ref{zac}) guarantees that equation (\ref{Aw1}) is satisfied in the strong sense (\ref{A1}).
Equation (\ref{A11}) is obtained by multiplying (\ref{A1}) by $b'(\vr)$.
To pass to the limit from inequality (\ref{galen})$_{(\vr_N,\vu_N)}$ to inequality (\ref{Aw3}) one uses, at the left-hand side the lower weak semi-continuity of convex functionals. The details of this limit passage are available in \cite[Section 4.3.4]{JiNoSIMA}.  We have thus established Lemma \ref{LLemmaA}.
}

\section{Relative energy for the approximate system}\label{5}
In this Section, we shall derive a convenient form of the relative energy inequality for any generalized solution of the approximate system (\ref{A1}--\ref{A4}).
This result is subject of the following lemma.
\bLemma{Renergy} Let all assumptions of Lemma \ref{LLemmaA} be satisfied. Let $(\vr=\vre,\vu=\vue)$ and an associted tensor field $\tn Z=\tn Z_\ep$ be a generalized solution to problem (\ref{A1}--\ref{A4})$_{\ep>0}$ constructed in Lemma \ref{LLemmaA}. Then $\vr$, $\vu$, $\tn Z$ satisfy the so called relative energy inequality
\begin{align}\label{rea}
&\intO{\Big(\frac 12\vr|\vc v-\vc V|^2
+E_\delta(\vr|r)\Big)(\tau,x)}
\nonumber\\
&\qquad+\int_0^\tau\int_{ \Gamma_{\rm in}} E_\delta(\vr_B|\vr)|\vu_B\cdot\vc n|{\rm d}S_x{\rm d}t
+ \int_0^\tau\int_{\Gamma_{\rm out}}E_\delta(\vr|r)\vu_B\cdot\vc n{\rm d}S_x{\rm d} t
\nonumber\\
&\qquad+\int_0^\tau\intO{\Big(\tn S(\Grad\vu):\Grad(\vc v-\vc V)
+\ep H_\delta''(\vr)|\Grad\vr|^2 { +\ep|\Grad\vc v|^4}\Big)}{\rm d}t
\nonumber\\
&\le
{ \intO{\Big(\frac 12\vr_0|\vc v_0-\vc V(0,\cdot)|^2+E_\delta(\vr_0|r(0,\cdot)\Big)}}
\nonumber\\
&\qquad+\int_0^\tau\int_{\Gamma_{\rm in}}\Big(H_\delta(r)-rH'_\delta(r) -H_\delta(\vr_B) +\vr_BH_\delta'(r)\Big)\vu_B\cdot\vc n{\rm d}S_x{\rm d} t
\nonumber\\
&\qquad+\int_0^\tau\int_\Omega\Big(\vr(\vc V-\vc v)\cdot\partial_t\vc V+\vr\vu\cdot\Grad\vc U\cdot(\vc V-\vc v)
\nonumber\\
&\qquad\qquad+\Big(p_\delta(r) -p_\delta'(r)(r-\vr)-p_\delta(\vr)\Big){\rm div}\vc U
\nonumber\\
&\qquad\qquad+\frac{r-\vr}r p'_\delta(r)(\vc v-\vc V)\cdot\Grad r - p'_\delta(r)\vc v\cdot\Grad r -p_\delta(r){\rm div}\vc V\Big){\rm d}x{\rm d}t
\nonumber\\
&\qquad+\int_0^\tau\intO{\Big(\tn Z:\Grad\vc V +\ep\Grad\vr\cdot\Grad(\vu- \vc V)\cdot\vc V\Big)
}{\rm d}t,
\end{align}
with a.e. $\tau\in \overline I$ and with any couple $(r,\vc U)$ belonging to class (\ref{Ugen}--\ref{rgen}).
In the above $\vc V=\vc U-\vu_\infty$, $\vc v=\vu-\vu_\infty$.
\eL
\bigskip\noindent
{\bf Proof of Lemma \ref{LRenergy}}\\
We take any couple $(r,\vc U)$ (and $\vc V$ related to $\vc U$) belonging to class (\ref{Ugen}--\ref{rgen}).
\begin{enumerate}
\item Using in the regularized continuity equation (\ref{Aw1}) the test function $\varphi=\frac 12|\vc V|^2$, we get:
\begin{align}\label{rea1}
&\intO{\frac12\vr(\tau,x)|\vc V|^2(\tau,x)} - \intO{\frac 12\vr_0(x)|\vc V|^2(0,x)}
\nonumber\\
&\qquad=\int_0^\tau\intO{\Big(\vr\vc V\cdot\partial_t\vc V + \vr \vu \cdot \Grad \vc V\cdot\vc V-\ep\Grad\vr\cdot\Grad\vc V\cdot\vc V\Big) }{\rm d}t
\end{align}
for all $\tau\in \overline I.$
\item Using in the momentum equation (\ref{Aw2}) the test function $\bfphi=-\vc V$, we obtain:
\begin{align}\label{rea2}
&-\intO{\vr\vc v(\tau,\cdot)\cdot\vc V(\tau,\cdot)} + \intO{\vr_0\vc v_0(\cdot)\vc V(0,\cdot)}
= { \int_0^\tau\intO{\tn Z_\ep:\Grad\varphi}{\rm d} t}
\nonumber\\
&
+\int_0^\tau\int_\Omega  \Big(-\vr\vc v\partial_t\vc V -\vr \vu \cdot \Grad \vc V\cdot \vu
{ +\vr}\vu\cdot\Grad(\vu_\infty\cdot \vc V)
{ - p_\delta(\vr) \Div \vc V}
\nonumber\\
&
+
\ep\Grad\vr\cdot\Grad\vu\cdot\vc V { -\ep\Grad\vr \cdot\Grad(\vu_\infty\cdot \vc V)}+\mathbb{S}(\Grad \vu) : \Grad \vc V\Big) {\rm d} x{\rm d}t
\end{align}
for all $\tau$ in $\overline I.$
\item Employing in the regularized continuity equation (\ref{Aw1}) the test function $\varphi=-H'(r)$, we get:
\begin{align}\label{rea3}
&-\intO{(\vr H_\delta'(r))(\tau,x)} + \intO{\vr_0(x) H_\delta'(r(0,x)}
\nonumber\\
&\qquad=\int_0^\tau\intO{\Big(-\vr \frac{p_\delta'(r)}r(\partial_t r +  \vu \cdot \Grad r)+\ep H_\delta''(r)\Grad\vr\cdot\Grad r\Big) }{\rm d}t
\nonumber\\
&\qquad+ \int_0^\tau\int_{\Gamma_{\rm in}} \vr_B H_\delta'(r) \vu_B \cdot \vc{n}   {\rm d}S_x{\rm d} t
+ \int_0^\tau\int_{\Gamma_{\rm out}} \vr H_\delta'(r) \vu_B\cdot \vc{n} \varphi \ {\rm d}S_x{\rm d} t
\end{align}
with any $\tau\in \overline I$, where we have used the identity (\ref{H+}) in the form
\begin{equation}\label{Hfd+}
r H'_\delta(r) - H(r)=p_\delta(r),
\end{equation}
cf. (\ref{pdelta}), (\ref{Hfd}).
\item The identity (\ref{Hfd+})  yields
\begin{equation}\label{rea4}
\intO{\Big(r H_\delta'(r)-H(r)\Big)}-\intO{\Big(r H_\delta'(r)-H_\delta(r)\Big)(0,\cdot)}=\int_0^\tau\intO{\partial_t p_\delta(r)}{\rm d}t
\end{equation}
for all $\tau\in \overline I$.
\item
Now, summing equations (\ref{rea1}), (\ref{rea2}), (\ref{rea3}), (\ref{rea4}) with the energy inequality (\ref{Aw3}) we get:
\begin{align}\label{rea-}
&\intO{\Big(\frac 12\vr|\vc v-\vc V|^2
+E_\delta(\vr|r)\Big)(\tau,x)}
\nonumber\\
&\qquad+\int_0^\tau\int_{ \Gamma_{\rm in}} E_\delta(\vr_B|\vr)|\vu_B\cdot\vc n|{\rm d}S_x{\rm d}t+
\int_0^\tau\int_{ \Gamma_{\rm out}}\Big( H_\delta(\vr)- \vr H_\delta'(r)  \Big)\vu_B\cdot\vc n{\rm d}S_x{\rm d}t
\nonumber\\
&\qquad+\int_0^\tau\intO{\Big(\tn S(\Grad\vu):\Grad(\vc v-\vc V)
+\ep H_\delta''(\vr)|\Grad\vr|^2 { +\ep|\Grad\vc v|^4}\Big)}{\rm d}t
\nonumber\\
&\le
{ \intO{\Big(\frac 12\vr_0|\vu_0-\vc U(0,\cdot)|^2+E_\delta(\vr_0|r(0,\cdot)\Big)}}
\nonumber\\
&\qquad+\int_0^\tau\int_{\Gamma_{\rm in}}\Big( -H_\delta(\vr_B) +\vr_BH_\delta'(r)\Big)\vu_B\cdot\vc n{\rm d}S_x{\rm d} t
\nonumber\\
&\qquad+\int_0^\tau\int_\Omega\Big(\vr(\vc V-\vc v)\cdot\partial_t\vc V+\vr\vu\cdot\Grad\vc U\cdot(\vc V-\vc v) +\Big(p_\delta(r)-p_\delta(\vr)\Big){\rm div}\vc U
\nonumber\\
&\qquad\qquad+\frac{r-\vr}r p'_\delta(r)\Big(\partial_t r+\vc u\cdot\Grad r\Big) - p'_\delta(r)\vc u\cdot\Grad r -p_\delta (r){\rm div}\vc U\Big){\rm d}x{\rm d}t
\nonumber\\
&\qquad+\int_0^\tau\intO{\Big(\tn Z_\ep:\Grad\vc V +\ep\Grad\vr\cdot\Grad(\vu- \vc V)\cdot\vc V
}{\rm d}t.
\end{align}
\item Now, we use in the before last line three elementary facts:
\begin{enumerate}
\item First, in view of (\ref{rgen}),
$$
\partial_t r+\vc u\cdot\Grad r= -r{\rm div}\vc U +(\vc u-\vc U)\cdot\Grad r
$$
\item Second, by definition of $\vc v$ and $\vc V$,
$$
- \vc u\cdot\Grad r=-\vc v\cdot\Grad r-\vc u_\infty\cdot\Grad r,\;-{\rm div}\vc U=-{\rm div}\vc V-{\rm div}\vc u_\infty.
$$
\item  Third, by virtue of Stokes formula,
$$
-\intO{\Big(p_\delta'(r)\vu_\infty\cdot \Grad r +p_\delta(r){\rm div}\vu_\infty\Big)}=-\intO{p_\delta(r)\vu_B\cdot\vc n}.
$$
\item
Consequently, the before last line reads
\begin{align*}
&\int_0^\tau\int_\Omega\Big(\frac{r-\vr}r p'_\delta(r)\Big(\partial_t r+\vc u\cdot\Grad r\Big) - p'_\delta(r)\vc u\cdot\Grad r -p_\delta(r){\rm div}\vc U\Big){\rm d}x{\rm d}t
\\
&
=\int_0^\tau\int_\Omega\Big((\vr-r)p'_\delta(r){\rm div} \vc U +\frac{r-\vr}r p'_\delta(r)(\vu-\vc U)\cdot\Grad r
{ -p_\delta'(r){\vc v}\cdot\Grad r}-p_\delta(r){\rm div}\vc V\Big){\rm d}x{\rm d}t
\\
&\qquad
-\int_{\Gamma_{\rm in}} p_\delta(r)\vc u_B\cdot\vc n{\rm d}S -\int_{\Gamma_{\rm out}} p_\delta(r)\vc u_B\cdot\vc n {\rm d}S.
\end{align*}
\end{enumerate}
\item Inserting the latter formula into (\ref{rea-}) and using also (\ref{Hfd}-\ref{Ed}), (\ref{Hfd+}), we obtain (\ref{rea}). Lemma \ref{LRenergy} is thus proved.
\end{enumerate}


\section{Limit $\ep\to 0$. Proof of Theorem \ref{TM1}: start}
\label{D}
The aim in this section is to pass to the limit in the weak formulation (\ref{Aw1}--\ref{Aw3}) of the problem (\ref{A1}--\ref{A4})$_{(\vr_\ep,\vu_\ep)}$ in order to
recover the weak formulation of problem (\ref{i1}--\ref{i5}) - cf. (\ref{fs}--\ref{me})$_{p=p_\delta}$, (\ref{ei})$_{p=p_\delta,H=H_\delta}$- and in the
relative energy inequality (\ref{rea})$_{(\vre,\vue,\tn Z_\ep)}$  to recover (\ref{rei})$_{p=p_\delta,H=H_\delta, E=E_\delta}$.

Estimates (\ref{ts1?}--\ref{ts2?+}) yield  uniform bounds
\begin{equation}\label{ts1}
\|\vr_\epsilon|\vu_\epsilon|^2\|_{L^\infty(I,L^1(\Omega ))}\le
 L({\rm data}),
\end{equation}
\begin{equation}\label{ts2}
\|\vu_\epsilon\|_{L^2(I,W^{1,2}(\Omega ))}\le L({\rm data}),
\end{equation}
\begin{equation}\label{ts3}
\|\vr_\epsilon\|_{L^\infty(I,L^\gamma(\Omega))}\le
L({\rm data}).
\end{equation}
Further, there is an $\delta$-dependent bound
\begin{equation}\label{ts4}
\delta^{1/\beta}\|\vr_\epsilon\|_{L^\infty(I,L^\beta(\Omega))}\le
L({\rm data})
\end{equation}
and $\ep$-dependent bounds
\begin{align}
&\ep^{1/4}\|\vu_\epsilon-\vu_\infty\|_{L^4(I,W^{1,4}(\Omega ))}\le L({\rm data},\delta),\;\|\tn Z_\ep\|_{L^{4/3}(Q_T)}\to 0,
\label{ts1plus}\\
&\sqrt \epsilon \|\nabla\vr_\epsilon\|_{L^2( Q_T)}\le
L({\rm data}).\label{ts8}
\end{align}
Moreover, the momentum equation (\ref{Aw2}) provides a  refined bound for the pressure, which reads
\begin{equation}\label{ts7}
\|\vr_\epsilon\|_{L^{\beta+1} ((0,T)\times K)}\le L({\rm data},\delta,K),\;\;\mbox{with any compacts $K\subset\Omega$},
\end{equation}
cf. \cite[Lemma 5.2]{JiNoSIMA}.

From estimates (\ref{ts3}--\ref{ts4}) we basically know, that there is a couple $(\vr,\vu)\in L^\infty((0,T);L^\beta(\Omega))
\times L^2((0,T); W^{1,2}(\Omega))$ which is a weak limit of a conveniently chosen subsequence of the sequence $(\vr_\ep,\vu_\ep)$
(not relabeled). It is proved in \cite{JiNoSIMA}, that this limit belongs to the class (\ref{fsw}) and satisfies the continuity equation
(\ref{ce}), the momentum equation  (\ref{me})$_{p=p_\delta}$ and the energy inequality (\ref{ei})$_{p=p_\delta, H=H_\delta}$. This nontrivial result
is obtained  through the key idea  asserting that
\begin{equation}\label{a.e.}
\vr_\ep\to\vr \;\mbox{a.e. in $Q_T$}\footnote{This is the key point in the existence theory of compressible Navier-Stokes equations, whatever is
the geometrical setting and whatever are the boundary conditions, cf. Lions \cite{LI4}, \cite{FNP}.}
\end{equation}

It remains to pass to the limit in the approximate relative energy inequality (\ref{rea})$_{(\vre,\vue,\tn Z_\ep)}$
and get (\ref{rei})$_{p=p_\delta,E=E_\delta, H=H_\delta}$. To this end we take $0<s<t<T$ and integrate  (\ref{rea})$_{(\vre,\vue,\tn Z_\ep)}$
over $\tau$ from $s$ to $t$.
We get
\begin{align}\label{rea+}
&\int_s^t\intO{\Big(\frac 12\vre|\vc v_\ep-\vc V|^2
+E_\delta(\vre|r)\Big)(\tau,x)}{\rm d}\tau
\nonumber\\
&
\qquad+\int_s^t\int_0^\tau\int_{ \Gamma_{\rm in}} E_\delta(\vr_B|\vre)|\vu_B\cdot\vc n|{\rm d}S_x{\rm d}t{\rm d}\tau
+ \int_s^t\int_0^\tau\int_{\Gamma_{\rm out}}E_\delta(\vre|r)\vu_B\cdot\vc n{\rm d}S_x{\rm d} t {\rm d}\tau
\nonumber\\
&\qquad+\int_s^t\int_0^\tau\intO{\Big(\tn S(\Grad\vue):\Grad(\vc v_\ep-\vc V)
+\ep H_\delta''(\vre)|\Grad\vre|^2 { +\ep|\Grad\vc v_\ep|^4}\Big)}{\rm d}t {\rm d}\tau
\nonumber\\
&\le
{(t-s) \intO{\Big(\frac 12\vr_0|\vc v_0-\vc V(0,\cdot)|^2+E_\delta(\vr_0|r(0,\cdot)\Big)}}
\nonumber\\
&\qquad+(t-s)\int_0^\tau\int_{\Gamma_{\rm in}}\Big(H_\delta(r)-rH'_\delta(r) -H_\delta(\vr_B) +\vr_BH_\delta'(r)\Big)\vu_B\cdot\vc n{\rm d}S_x{\rm d} t {\rm d}\tau
\nonumber\\
&\qquad+\int_t^s\int_0^\tau\int_\Omega\Big(\vre(\vc V-\vc v_\ep)\cdot\partial_t\vc V+\vre\vue\cdot\Grad\vc U\cdot(\vc V-\vc v_\ep)
\nonumber\\
&\qquad\qquad+\Big(p_\delta(r) -p_\delta'(r)(r-\vre)-p_\delta(\vre)\Big){\rm div}\vc U
\nonumber\\
&\qquad\qquad+\frac{r-\vr}r p'_\delta(r)(\vc v_\ep-\vc V)\cdot\Grad r - p'_\delta(r)\vc v_\ep\cdot\Grad r -p_\delta(r){\rm div}\vc V\Big){\rm d}x{\rm d}t{\rm d}\tau
\nonumber\\
&\qquad+\int_t^s\int_0^\tau\intO{\Big(\tn Z_\ep:\Grad\vc V +\ep\Grad\vre\cdot\Grad(\vue- \vc V)\cdot\vc V\Big)
}{\rm d}t{\rm d}\tau.
\end{align}

The relation (\ref{a.e.}) in combination with estimates (\ref{ts1}--\ref{ts7}) employed in (\ref{rea+})
allows to pass to the limit in all of its terms except the term
\begin{equation}\label{reap}
-\int_t^s\int_0^\tau\intO{p_\delta(\vre){\rm div}\vc U}{\rm d}t{\rm d}\tau
\end{equation}
 Indeed:
\begin{enumerate}
\item  Under the assumption $p'(\vr)>0$ the Helmholtz function $H_\delta$ is strictly convex.
We can therefore omit  at the left hand side of (\ref{rea}) all integrals over the boundary, since they are non negative.
The same is true for the right hand side volumic integrals containing $\ep$ as a multiplier.
 \item We use lower weak semicontinuity of the convex functionals together with standard compactness arguments\footnote{Here and hereafter, the standard compactness arguments include the Sobolev imbeddings, Arzela-Ascoli theorem and Lions-Aubin lemma.}  in  the remaining terms
of the left hand side of (\ref{rea+}).
\item We employ estimate (\ref{ts8}) to get rid of the term containing multiplication by $\ep$ at the right hand side,
and relation (\ref{ts1plus}) to get rid of the term containing $\tn Z_\ep$.
\item We use the standard compactness arguments  in all terms of the right hand side except the term (\ref{reap}). The term (\ref{reap})
cannot be so far treated, due to the fact that estimate (\ref{ts7}) does not hold up to the boundary.
\end{enumerate}

To remedy to the above problem, we intend to prove that the integral of the pressure over the space-time cylinder whose basis is an inner neighborhood of the boundary $\partial\Omega$, is comparable with the measure of the neighborhood in a certain positive power. To this end, we show first the following lemma.

\begin{Lemma}\label{localp}
Let $\Omega$ be a bounded Lipschitz domain, $\hat U_h^-$, $h>0$ an inner neighborhood of its boundary - see (\ref{hatU}) - and $\alpha,\kappa>1$. Consider
a sequence
$
(p_\ep,\vc z_\ep,\vc F_\ep,\tn G_\ep)_{\ep>0}
$
of functions  which satisfy equation
\begin{equation}\label{distr}
\partial_t\vc z_\ep+\vc F_\ep+{\rm div}\tn G_\ep +\Grad p_\ep=0\;\mbox{in ${\cal D}'(Q_T;\R^d)$}.
\end{equation}
Suppose finally that $p_\ep\in L^1(Q_T)$, while
$$
(z_\ep,\vc F_\ep,\tn G_\ep)_{\ep>0}\; \mbox{is bounded in $L^\infty(0,T;L^\alpha(\Omega))\times L^\kappa(Q_T)\times L^\kappa(Q_T)$ by $k>0$}
$$
uniformly with respect to $\ep$.

Then there exists $h_0>0$ and $c=c(k,T,\Omega)>0$ such that
\begin{equation}\label{localpe}
\int_0^T\int_{\hat U_h^-}p(\vre){\rm d}x{\rm d}t\le ch^{\Gamma},\;\mbox{where $\Gamma=\min\{1/\alpha', 1/\kappa'\}$},
\end{equation}
for all $0<h<h_0$ uniformly with respect to $\ep$.
\end{Lemma}
\noindent
{\bf Proof of Lemma \ref{localp}}\\
We shall proceed in several steps.
\begin{enumerate}
\item
We shall start the proof by recalling the so called Bogovskii lemma, see e.g. Galdi \cite{Galdi} or \cite{NOST}.
\bLemma{Bog}
Let $\Omega$ be a bounded Lipschitz domain. Then there exists a linear operator
$$
{\cal B}:\{f\in C^\infty_c(\Omega;\R^3)\,|\,\intO{f}=0\}\mapsto C^\infty_c(\Omega;\R^3)\}
$$
such that:
\begin{enumerate}
\item
${\rm div}{\cal B}[f]=f$ \item
${\cal B}$ is bounded linear operator from $\overline{L}^p(\Omega)$ to $W^{1,p}(\Omega)$
for any $1<p<\infty$ (i.e. there is $c=c(p)>0$ such that $\|{\cal B}[f]\|_{W^{1,p}(\Omega;\R^3)}\le c(p)\|f\|_{L^p(\Omega)}$ for all
$f\in \overline L^p(\Omega)$).
In the above $\overline L^p(\Omega):=\{f\in L^p(\Omega)\,|\,\intO{f}=0\}$.
\end{enumerate}
\eL
\item Let $h_0$ be sufficiently small, such that $\hat U_h^-\subset\Omega$ and $|\hat U_h^-|<|\Omega|/2$.
We take in Lemma \ref{LBog}
$$
f(x)=1_{\Omega\setminus\hat U_h^-} (x)-\frac{|\Omega|-|\hat U_h^-|}{|\Omega|},
$$
 where $1_A$ denotes the characteristic
function of the set $A$. We observe that: 1) Since $\Omega$ is Lipschitz, there exists $\overline c>0$ such that for all $0<h<h_0$,
\begin{equation}\label{Uh}
|\hat U_h^-|\le\overline c h.
\end{equation}
2) Function $f$ satisfies all assumptions of Lemma \ref{LBog} with any $1<p<\infty$. Consequently, setting $\psi(x)={\cal B}[f]$, we get estimate
\begin{equation}\label{Boge}
\|\psi\|_{W^{1,p}(\Omega)}\le c(p) h^{1/p}\; \mbox{for all $0<h<h_0$}.
\end{equation}
\item Now, we use in equation (\ref{distr}) as the test function the function
$$
\varphi(t,x)=\eta(t)\psi(x),\;\mbox{where $0\le\eta\in C^1_c((0,T))$, $\eta(t)=1$ in $(\Delta, T-\Delta)$, $\Delta\in (0,T/2)$, $|\partial_t\eta|\le 2/\Delta$}.
$$
(This is an admissible test function, as one can show by a density argument.)
We get
\begin{equation}\label{konec}
-\int_0^T\eta\intO{p{\rm div}\varphi}{\rm d}t=I_1+I_2+I_3,
\end{equation}
where
$$
I_1=\int_0^T\partial_t\eta\intO{\vc z_\ep\cdot\psi}{\rm d}t,\; I_2=\int_0^T\eta\intO{\tn G_\ep:\Grad\psi}{\rm d}t,\;I_3=-\int_0^T\eta\intO{\vc F_\ep\cdot\psi}{\rm d}t.
$$
Seeing that
$$
-\int_0^T\eta\intO{p{\rm div}\varphi}{\rm d}t=\frac{|\Omega|-|\hat U_h^-|}{|\Omega|}\int_0^T\eta\int_{\hat U_h^-}p(\vre){\rm d} x{\rm d}t
-\frac{|\hat U_h^-|}{|\Omega|}\int_{\Omega\setminus \hat U_h^-}p(\vre){\rm d} x{\rm d}t
$$
and that
$$
|I_1|\le c\|\vc z_\ep\|_{L^\infty(0,T;L^\alpha(\Omega))}\|\psi\|_{L^{\alpha'}(\Omega)},\;|I_2|+|I_3|\le c
\|\vc F_\ep,\tn G_\ep\|_{L^\kappa(Q_T)} \|\psi\|_{L^{\kappa'}(\Omega)}
$$
with $c$ independent of $\eta$, i.e., in particular, independent of $\Delta$. Inserting both latter observations to (\ref{konec}) while taking into account
(\ref{Uh}) and (\ref{Boge}) completes the proof of Lemma \ref{localp}.
\end{enumerate}

We are now ready to treat the problematic term (\ref{reap}).
We apply Lemma \ref{localp} to the momentum equation (\ref{Aw2}), i.e., we set $p_\ep=p_\delta(\vr_\ep)$, $\vc z_\ep=\vr_\ep\vc v_\ep$,
$\vc F_\ep=\vre\vue\cdot\Grad\vu_\infty+\ep\Grad\vre\cdot\Grad\vc v_\ep$, $\tn G_\ep=-\tn Z_\ep+\vre\vue\otimes\vc v_\ep-\tn S(\Grad\vue)$.
Due to (\ref{ts1}), (\ref{ts4}), $\vc z_\ep$ is bounded in $L^\infty(0,T; L^{\frac{2\beta}{\beta+1}}(\Omega))$, by virtue of
(\ref{ts1}), (\ref{ts2}), (\ref{ts4}),(\ref{ts1plus}), $\tn G_\ep$ is bounded in $L^{\min\{4/3,\frac{6\beta}{4\beta+3}\}}(Q_T)$, and
finally, due to (\ref{ts1}), (\ref{ts4}), (\ref{ts1plus}) and (\ref{ts8}), $\vc F_\ep$ is bounded in $L^{4/3}(Q_T)$.
We thus obtain
\begin{equation}\label{herror}
\limsup_{\ep\to 0}|\int_t^s\int_0^\tau\int_{\hat U_h^-}{p_\delta(\vre){\rm div}\vc U}{\rm d}x{\rm d}t{\rm d}\tau|\le
c\sup_{(t,x)\in Q_T}|{\rm div}\vc U(t,x)|\; h^\Gamma\;\mbox{with some $\Gamma>0$},
\end{equation}
for all $0<h<h_0$ and $0\le t<s\le T$.

We may write with any $0<h<h_0$,
\begin{equation}\label{decomp}
-\int_t^s\int_0^\tau\intO{p_\delta(\vre){\rm div}\vc U}{\rm d}t{\rm d}\tau=
\int_t^s\int_0^\tau\int_{\hat U_h^-}{p_\delta(\vre){\rm div}\vc U}{\rm d}x{\rm d}t{\rm d}\tau+
\int_t^s\int_0^\tau\int_{\Omega\setminus\hat U_h^-}{p_\delta(\vre){\rm div}\vc U}{\rm d}x{\rm d}t{\rm d}\tau,
\end{equation}
where due to (\ref{ts7}) and (\ref{a.e.}),
$$
\lim_{\ep\to 0}
\int_t^s\int_0^\tau\int_{\Omega\setminus\hat U_h^-}{p_\delta(\vre){\rm div}\vc U}{\rm d}x{\rm d}t{\rm d}\tau
=
\int_t^s\int_0^\tau\int_{\Omega\setminus\hat U_h^-}{p_\delta(\vr){\rm div}\vc U}{\rm d}x{\rm d}t{\rm d}\tau
$$
with any $0<h<h_0$.
We already know that, in particular, $p(\vr)\in L^1(Q_T)$. Therefore
$$
\lim_{h\to 0}\int_t^s\int_0^\tau\int_{\Omega\setminus\hat U_h^-}{p_\delta(\vr){\rm div}\vc U}{\rm d}x{\rm d}t{\rm d}\tau=
\int_t^s\int_0^\tau\int_{\Omega}{p_\delta(\vr){\rm div}\vc U}{\rm d}x{\rm d}t{\rm d}\tau.
$$
Using the both latter facts and estimate (\ref{herror}) in the decomposition (\ref{decomp}) we get the desired conclusion,
namely, that
$$
\int_t^s\int_0^\tau\intO{p_\delta(\vre){\rm div}\vc U}{\rm d}t{\rm d}\tau\to \int_t^s\int_0^\tau\intO{p_\delta(\vr){\rm div}\vc U}{\rm d}t{\rm d}\tau
$$
as $\ep\to 0$. This is the last element needed to get from inequality (\ref{rea+}),
\begin{align}\label{rea++}
&\int_s^t\intO{\Big(\frac 12\vr|\vc v-\vc V|^2
+E_\delta(\vr|r)\Big)(\tau,x)}{\rm d}\tau
\nonumber\\
&\le
{(t-s) \intO{\Big(\frac 12\vr_0|\vc v_0-\vc V(0,\cdot)|^2+E_\delta(\vr_0|r(0,\cdot)\Big)}}
\nonumber\\
&\qquad+(t-s)\int_0^\tau\int_{\Gamma_{\rm in}}\Big(H_\delta(r)-rH'_\delta(r) -H_\delta(\vr_B) +\vr_BH_\delta'(r)\Big)\vu_B\cdot\vc n{\rm d}S_x{\rm d} t {\rm d}\tau
\nonumber\\
&\qquad+\int_t^s\int_0^\tau\int_\Omega\Big(\vr(\vc V-\vc v)\cdot\partial_t\vc V+\vr\vu\cdot\Grad\vc U\cdot(\vc V-\vc v)
\nonumber\\
&\qquad\qquad+\Big(p_\delta(r) -p_\delta'(r)(r-\vr)-p_\delta(\vr)\Big){\rm div}\vc U
\nonumber\\
&\qquad\qquad+\frac{r-\vr}r p'_\delta(r)(\vc v-\vc V)\cdot\Grad r - p'_\delta(r)\vc v\cdot\Grad r -p_\delta(r){\rm div}\vc V\Big){\rm d}x{\rm d}t{\rm d}\tau
\end{align}
for all $0\le s<t\le T$. The final step is to devide (\ref{rea++}) by $t-s$ and effectuate the limit $t\to s$. The theorem on Lebesgue points
then guarantees the satisfaction of the relative energy inequality (\ref{rei}) with $(p,H,E)$ replaced by $(p_\delta, H_\delta, E_\delta)$.

\section{Limit $\delta\to 0$. Proof of Theorem \ref{TM1}: end}\label{sE}
Our final goal is to pass to the limit in the weak formulation (\ref{i1}--\ref{i5})$_{p=p_\delta,\vr_\delta,\vu_\delta}$-
cf. (\ref{fs}--\ref{me})$_{p=p_\delta,\vrd,\vud}$, (\ref{ei})$_{p=p_\delta,H=H_\delta,\vrd,\vud}$, in order to
recover the weak formulation (\ref{fs}--\ref{me}), (\ref{ei}) of problem (\ref{i1}--\ref{i5}) - and in the
relative energy inequality (\ref{rei})$_{p=p_\delta, H=H_\delta, E=E_\delta,\vrd,\vud}$ in order  to recover (\ref{rei}).

Estimates (\ref{ts1}--\ref{ts4}) yield uniform bounds
\begin{equation}\label{ts1+}
\|\vr_\delta\vud^2\|_{L^\infty(I,L^1(\Omega))}\le
L({\rm data}),
\end{equation}
\begin{equation}\label{ts3+}
\|\vu_\delta\|_{L^2(I,W^{1,2}(\Omega;R^3 ))}\le L({\rm data}),
\end{equation}
and $\delta$ dependent bound
 \begin{equation}\label{ts4+}
\delta^{1/\beta}\|\vr_\delta\|_{L^\infty(I,L^\beta(\Omega))}\le
L({\rm data}).
\end{equation}
Moreover, the momentum equation provides, under condition $\gamma>3/2$, a  refined bound for pressure, which reads
\begin{equation}\label{ts7+}
\|\vr_\delta\|_{L^{\gamma+\alpha(\gamma)} ((0,T)\times K)}\le L({\rm data},K),\;\;\mbox{with any compacts $K\subset\Omega$},
\end{equation}
and with $\alpha=\min\{\frac 23\gamma-1,\frac\gamma 2\}$, cf. \cite[Lemma 6.2]{JiNoSIMA}.

From estimates (\ref{ts3}--\ref{ts4}) we deduce, that there is a couple $(\vr,\vu)\in L^\infty((0,T);L^\gamma(\Omega))
\times L^2((0,T); W^{1,2}(\Omega))$
which is a weak limit of a conveniently chosen subsequence of the sequence $(\vrd,\vud)$
(not relabeled).  We also know that
\begin{equation}\label{strd}
\vrd\to\vr \;\mbox{a.e. in $Q_T$}
\end{equation}
provided $\gamma>3/2$. As underlined in the previous section, the latter convergence relation is crucial in the theory of compressible
Navier-Stokes equations and its proof is quite involved. We refer to \cite[Section 6]{JiNoSIMA} for the proof in the present context.
It is shown, that this limit belongs to the class (\ref{fs}) and satisfies the continuity equation
(\ref{ce}), the  momentum equation  (\ref{me}) and the energy inequality (\ref{ei}).
Our task is only to pass to the limit
in the relative energy inequality (\ref{rei})$_{p=p_\delta, H=H_\delta,E=E_\delta,\vrd,\vud}$ and to obtain (\ref{rei}).

Reasoning as in the previous section we discover that the only problematic term is \\$\int_0^\tau\intO{p_\delta(\vrd){\rm div}\vc U}{\rm d}t$.
Indeed, since estimate (\ref{ts7+}) is only local, it is not clear $p_\delta(\vrd)$ converges weakly to $p_\delta(\vr)$ "near the boundary",
even if we know (\ref{strd}). We will treat this difficulty exactly in the same manner as we have treated the term (\ref{reap}) in the previous
section, by using Lemma \ref{localp}. This lemma will be applied to the momentum equation (\ref{me})$_{p=p_\delta,\vrd,\vud}$ yielding estimate
(\ref{herror}) with $\vre$ replaced by $\vrd$. In this calculation, again, the value $3/2$ is a threshold that cannot be achieved.

The rest of the reasoning is the same as in the previous section and thus left to the reader. Theorem \ref{TM1} is proved.

\section{Stability and weak-strong uniqueness: Proof of Theorem \ref{TWS1}} \label{WS}

In this section we shall prove Theorem \ref{TWS1}. We shall show that the strong solutions to the problem (\ref{i1}-\ref{i5}) are stable in the class of dissiptive weak solutions. In particular,
any  dissipative weak solution of the  problem (\ref{i1}-\ref{i5}) coincides with the strong solution of the same problem emanating from the same initial
data and the same boundary conditions.

\subsection{Relative energy inequality with a strong solution as a test function}

If the test functions $(r,\vc U)$ in the relative energy inequality (\ref{rei}) obey equations
(\ref{i1}-\ref{i2}) almost everywhere in $Q_T$ the right hand side of the relative energy
becomes quadratic in differences $(\vr-r,\vu-\vc U)$. This observation is subject of the following lemma:

\bLemma{WS1}
Let $\Omega$ be a bounded Lipschitz domain.  Suppose that the pressure satisfies assumptions (\ref{pressure}).\footnote{This Lemma holds without condition
$p(0)=0$, $p'(\vr)>0$; regularity assumption of (\ref{pressure}) is enough.}
Let $(\vr,\vc u)$ be a dissipative weak solution to the Navier-Stokes equations (\ref{i1}-\ref{i5})  emanating from the finite energy initial data $(\vr_0,\vu_0)$
in class (\ref{feid}) and
boundary data $(\vr_B,\vu_B)$ in class (\ref{M1}), corresponding to extension $\vu_\infty$ of $\vu_B$- cf. (\ref{Le}). Let $(r,\vc U)$
belonging to the class
\begin{align}\label{s1-}
&0<\underline r\le r\le\overline r<\infty;\quad\vc U \in L^\infty((0,T)\times\Omega)
\nonumber\\
&\partial_t r,\partial_t\vc U,\Grad r,\Grad\vc U\in L^2(0,T; C(\overline\Omega))
\end{align}
be a strong solution of the same equations with  initial data $(r(0),\vc U(0))=(r_0,\vc U_0)$ and boundary data $(r|_{\Gamma_{\rm in}},
\vc U|_{\partial\Omega})=(r_B,\vu_B)$.
Then the relative energy inequality (\ref{rei}) takes the form:
\begin{align}\label{reis}
&\intO{\Big(\frac 12\vr|\vc v-\vc V|^2+E(\vr|r)\Big)(\tau)}+\int_0^\tau\intO{\tn S(\Grad(\vc v -\vc V)):\Grad(\vc v-\vc V)}{\rm d}t
\nonumber\\
&\qquad\le
{ \intO{\Big(\frac 12\vr_0|\vc v_0-\vc V_0|^2+E(\vr_0|r_0)\Big)}} + {\cal R}(\vr,\vc v|r,\vc V),
\end{align}
for a.e. $\tau\in I$. In the above, the remainder reads
$$
{\cal R}(\vr,\vc v|r,\vc V)=
\int_0^\tau\int_{\Gamma_{\rm in}}\Big( H(r_B)-H(\vr_B)+(\vr_B-r_B)H'(r_B)\Big)\vu_B\cdot\vc n{\rm d}S{\rm d}t
$$
$$
+\int_0^\tau\intO{\Big(\vr-r)(\vc V-\vc v)\cdot\partial_t\vc U + (\vr-r)\vc U\cdot\Grad\vc U\cdot(\vc V-\vc v)+
\vr(\vc v-\vc V)\cdot\Grad\vc U\cdot(\vc V-\vc v)\Big)}{\rm d}t
$$
$$
+
\int_0^\tau\intO{\Big(p(r)
 -p'(r)(r-\vr)-p(\vr)\Big){\rm div}\vc U }{\rm d}t
+\int_0^\tau\intO{\Big((1-\frac\vr r)p'(r)(\vc v-\vc V)\cdot\Grad r\Big)}{\rm d}t
$$
and
$\vc v=\vu-\vu_\infty$, $\vc v_0=\vu_0-\vu_\infty$, $\vc V=\vc  U-\vc u_\infty$, $\vc V_0=\vc U_0-\vc u_\infty$.
\eL
 \noindent
{\bf Proof of Lemma \ref{LWS1}}\\
We start by the observation that due to the regularity (\ref{s1-}) the couple $(r,\vc V)$  satisfies
\begin{equation}\label{i4+}
\partial_t r + {\rm div} (r \vc U) = 0\;\mbox{a.e. in $(0,T)\times\Omega$},
\end{equation}
\begin{equation}\label{i5++}
r\partial_t  \vc U + r\vc U\cdot\nabla\vc U + \nabla p(r) = {\rm div}{\tn S}(\nabla\vc U)\;\mbox{ a.e. in $(0,T)\times\Omega$}.
\end{equation}
{ Multiplying (\ref{i5++}) scalarly by $\vc v-\vc V$ and integrating over $\Omega$, we get
\begin{equation}\label{re1} \intO{\Big(r \partial_{t} \vc U + r \vc U \cdot \nabla \vc U + \nabla p(r)
\Big)\cdot(\vc v-\vc V)}{\rm d} t+ \intO{ \tn{S}(\nabla \vc U) : \nabla( \vc v - \vc V) }{\rm d}t = 0,
\end{equation}
where we have used the integration by parts in the last integral.

Next we calculate
$$
\intO{\Grad p(r)\cdot(\vc v-\vc V)}=\intO{p'(r)\vc v\cdot\Grad r}
+\intO{p(r){\rm div}\vc V}.
$$

Adding (\ref{i5++}) to the inequality (\ref{rei}) while taking into acount the above identity and relations (\ref{Hfd})
between $p$ and $H$, we obtain the inequality (\ref{reis}). Lemma \ref{LWS1} is proved.

\subsection{Two algebraic relations}

It is evident that  under the assumption $p'(\vr)>0$, the Helmholtz function $H$ is strictly convex and therefore
$E(\vr|r)>0$ for all $\vr\ge 0$, $r>0$, and $E(\vr|r)=0$ if and only if $\vr=r$. We can, however, prove more:

\bLemma{relaxed2}
Let $0<a<b<\infty$ and let $p$ satisfies (\ref{pressure}).
Then there exists a number $c=c(a,b)>0$ such that for all $\vr\in [0,\infty)$ and $r\in [a,b]$,
\begin{equation}\label{E1}
E(\vr|r)\ge c(a,b)\Big( 1_{{\cal O}_{\rm res}}(\vr) + \vr 1_{{\cal O}_{\rm res}}(\vr)+ (\vr-r)^2 1_{{\cal O}_{\rm ess}}(\vr)\Big),
\end{equation}
where $E$ is defined in (\ref{E}) and
\begin{equation}\label{calO}
{\cal O}_{\rm ess}=[a/2, 2b],\; {\cal O}_{res}= [0,\infty)\setminus {\cal O}_{\rm ess}.
\end{equation}
\eL
{\bf Proof}\\

If $\vr\in [a/2, 2b]$ we use the strict convexity of $H$ to obtain that
$$
E(\vr|r)\ge c|\vr -r|^2\;\mbox{where $c= c(a,b)>0$}.
$$

If $\vr\in R_+\setminus[a/2,2b]$, we observe that
$$
\partial_\vr E(\vr|r)= { H}'(\vr)-{ H}'(r),\quad\partial_r E(\vr|r)={ H}''(r)(r-\vr),
$$
where $s\to{ H}'(s)$ is an increasing function on $(0,\infty)$.
Now, relying on the monotonicity of functions $s\to E(s|r)$ and $s\to E(\vr|s)$ induced by the above formulas, we consider two situations. 1) If $\vr>
2b$, we observe that $E(\vr|2b)>0$, whence
${ H}(\vr)+ p(2b)> { H}'(2b)\vr$.
Consequently,
$$
p(2b)-p(b)+E(\vr|r)\ge p(2b)-p(b)+ E(\vr|b)= H(\vr)+ p(2b)-H'(b)\vr
\ge \Big(H'(2b)-H'(b)\Big)\vr.
$$
This inequality and the fact that $E(\vr,r)\ge E(2b,b)>0$, $p(2b)>p(b)$, $H'(2b)>H'(b)$ yield
$$
E(\vr|r)\ge c(1+\vr)
$$
with some $c=c(b)>0$.
%
%
2) If $\vr <a/2$ then
$$
E(\vr|r)\ge E(a/2| a) \ge \frac { E(a/2| a)} a \vr + \frac { E(a/2| a)} 2\ge c(1+\vr)
$$
with some $c=c(a)>0$. Lemma \ref{Lrelaxed2} is proved.
\bigskip

Since $\vr\mapsto p(\vr)$ is bounded on any compact subset of $[0,\infty)$ and since according to the definition of $E(\cdot|\cdot)$,
$H(\vr)\le c \Big(E(\vr| r) + 1 +\vr\Big)$ with any $\vr\in [0,\infty)$ and any $r\in [a/2,2b]$, where $c=c(a,b)>0$, we deduce
from estimate (\ref{E1}) the following result:
\bCorollary{relaxed3}
Suppose that $p$ assumption (\ref{pressure}) and (\ref{w1+}). Then for any $\vr\in [0,\infty)$ and $r\in [a/2, 2b]$ we have,
in addition to inequality (\ref{E1}),
\begin{equation}\label{E2}
p(\vr) 1_{{\cal O}_{res}}(\vr)\le c E(\vr|r)
\end{equation}
with some $c=c(a,b)>0$.
\eC

\subsection{Proof of Theorem \ref{TWS1}}
\subsubsection{The Gronwall inequality}\label{mainideas}

 The goal now is
 to find an estimate of the left hand side of
(\ref{reis})
from below by
\begin{equation}\label{how1+}
{ c\int_0^\tau\|{\vc v}-\vc V\|^2_{W^{1,2}(\Omega;\R^3)}{\rm d} t- \overline c'\int_0^\tau{\cal E}(\vr, \vc v\Big|r,\vc V){\rm d}t +{\cal E}(\vr, \vc v\Big|r,\vc V)\Big|_0^\tau,
}
\end{equation}
and the right hand side from above by
\begin{equation}\label{how2+}
\overline c\tau\|\vr_B-r_B\|_{L^1(\Gamma_{\rm in})}+
\delta\int_0^\tau\|{\vc v}-\vc V\|^2_{W^{1,2}(\Omega)}{\rm d}t
+ c'(\delta)\int_0^\tau a(t) {\cal E}(\vr, \vc v\Big|r,\vc V){\rm d}t
\end{equation}
with any $\delta>0$,
where $c,\overline c>0$ are independent of $\delta$, $\overline c'\ge 0$, $c'=c'(\delta)>0$, and $a\in L^1(0,T)$.
This process leads to the estimate
\begin{equation}\label{how3}
{\cal E}(\vr, \vc v\Big|r,\vc V)(\tau)\le {\cal E}(\vr_0, \vc v_0\Big|r(0),\vc V(0)) +\overline c T \|\vr_B-r_B\|_{L^1(\Gamma_{\rm in})}
+ c\int_0^\tau a(t) {\cal E}(\vr, \vc v\Big|r,\vc V){\rm d}t
\end{equation}
that implies estimate (\ref{stability1+}) by  the Gronwall inequality. In the rest of this section, we shall perform this program.

We start by observing that the bound from below (\ref{how1+}) holds true with $\overline c'=0$ and some $c=c(\mu)>0$. Indeed, to see this
we may use the Korn type inequality
$$\intO{\tn S(\nabla\vc w):\nabla\vc w}\ge \overline c\|\nabla\vc w\|^2_{L^2(\Omega;\R^3)}$$
holding for all $\vc w\in W_0^{1,2}(\Omega;\R^3)$ with $\overline c>$ independent of $\vc w$
(which can be, in this setting, easily proved by using the integration by parts and a density argument)
and the standard Poincar\'e inequality to deduce that
\bFormula{belowbound}
c\int_0^\tau\|{\vc v}-\vc V\|^2_{W^{1,2}(\Omega;\R^3)}{\rm d} t +{\cal E}(\vr, \vc v\Big|r,\vc V)\Big|_0^\tau
\le \int_0^\tau\intO{\tn S(\nabla(\vc v-\vc V):\nabla(\vc v-\vc V) }{\rm d }t+{\cal E}(\vr, \vc v\Big|r,\vc V)\Big|_0^\tau.
\eF

\subsection{Estimates of the remainder}
We introduce essential and residual sets in $\Omega$. To this end we take  in (\ref{calO}) $a=\underline r$, $b=\overline r$ and define for a.e. $t\in (0,T)$ the residual and essential subsets of $\Omega$
as follows:
\begin{equation}\label{essres}
N_{\rm ess}(t)=\{x\in\Omega\,\Big| \vr(t) \in {\cal O}_{\rm ess}\},\;
N_{\rm res}(t)= \Omega\setminus N_{\rm ess}(t).
\end{equation}

With this definition at hand and having assumption (\ref{w1+}) in mind, we deduce from Lemma \ref{Lrelaxed2} and Corollary \ref{Crelaxed3}
\begin{equation}\label{rentropy}
c\intO{\Big(\Big[1\Big]_{\rm res}+\Big[\vr\Big]_{\rm res}+\Big[p(\vr)\Big]_{\rm res}+\Big[\vr-r\Big]^2_{\rm ess}\Big)}\le
\intO{{ E}(\vr,\vu\Big| r,\vc U)}
\end{equation}
with some $c=c(\underline r,\overline r)>0$, where
we denote, for a function $h$ defined a.e. in $(0,T)\times\Omega$,
$$
[h]_{\rm ess}= h 1_{N_{\rm ess}},\; [h]_{\rm res}= h 1_{N_{\rm res}}.
$$

We are now in position to estimate the remainder ${\cal R}$ at  right hand side of the relative energy inequality (\ref{reis}). We shall do it in five steps.\\ \\
{\bf Step 1:}{\it The surface integral in ${\cal R}$}\\
We have immediately by the Taylor formula
\begin{equation}\label{rem0}
\Big|\int_0^\tau\int_{\Gamma_{\rm in}}\Big( H(r_B)-H(\vr_B)+(\vr_B-r_B)H'(r_B)\Big)\vu_B\cdot\vc n\Big){\rm d}S{\rm d}t\Big|\le
T\,c(|H'|_{C([\underline{\mathfrak{r}},\overline{\mathfrak{r}}])},|\vu_B|_{C(\partial\Omega)})\,\|\vr_B-r_B\|_{L^1(\Gamma_{\rm in})}.
\end{equation}
{\bf Step 2:} {\it The first volume integral in ${\cal R}$ }\\
We shall first estimate the ``essential part'' of the first two terms:
\begin{align}\label{remainder1}
&\int_0^\tau\intO{ [1]_{\rm ess}(\vr -r)(\partial_t \vc U + \vc U \cdot \Grad \vc U) \cdot (\vc V - \vc v)}{\rm d}t
\nonumber\\
&\qquad\le
\int_0^\tau\Big\|\partial_t\vc U+\vc U\cdot\Grad\vc U\|_{L^\infty(\Omega;\R^3)}\Big\|\Big[ \rho -r\Big]_{\rm ess}\Big\|_{L^2(\Omega)}
\Big\|{\vc v}-\vc V\Big\|_{L^2(\Omega;R^3)} {\rm d}t
\nonumber\\
&\qquad\le \delta \int_0^\tau \Big\|{\vc v}-\vc V\Big\|^2_{L^2(\Omega;\R^3)}{\rm d}t + c(\delta,\underline r,\overline r) \int_0^\tau a(t) {\cal E}\Big(\vr,\vc v\Big| r,\vc V\Big){\rm d}t,
\end{align}
where
$$
a =\|\partial_t\vc U+\vc U\cdot\Grad\vc U\|_{L^\infty(\Omega;\R^3)}^2\in L^1(0,T).
$$
Concerning the "residual part", we shall estimate the integrals over the sets $\{\vr\le\underline r/2\}$
and $\{\vr\ge 2\overline r\}$
separately.

\begin{align*}
&\int_0^\tau\intO{ 1_{\{\vr\le\underline r/2\}}(\rho -r)(\partial_t \vc U + \vc U \cdot \Grad \vc U) \cdot (\vc V - \vc v)}{\rm d}t
\nonumber\\
&\qquad\le 2\overline r \int_0^\tau\intO{ 1_{\rm res}\Big|\partial_t \vc U + \vc U \cdot \Grad \vc U\Big| \,\Big|\vc V - \vc v\Big|}{\rm d}t
\nonumber\\
&\qquad\le 2\overline r
\int_0^\tau\Big\|\partial_t\vc U+\vc U\cdot\Grad\vc U\|_{L^\infty(\Omega;\R^3)}\Big\|1_{\rm res}\Big\|_{L^2(\Omega)}
\Big\|{\vc v}-\vc V\Big\|_{L^2(\Omega;\R^3)} {\rm d}t
\nonumber\\
&\qquad\le \delta \int_0^\tau \Big\|{\vc v}-\vc V\Big\|^2_{L^2(\Omega;\R^3)}{\rm d}t + c(\delta,\underline r,\overline r) \int_0^\tau a(t) {\cal E}\Big(\vr,\vc v\Big| r,\vc V\Big){\rm d}t,
\end{align*}
where $a$ is given in (\ref{remainder1}).

Finally,
\begin{align*}
&\int_0^\tau\intO{ 1_{\{\vr\ge 2\overline r \}}(\vr)(\vr -r)(\partial_t \vc U + \vc U \cdot \Grad \vc U) \cdot (\vc V - \vc v)}{\rm d}t
\nonumber\\
&\qquad\le 2\int_0^\tau\intO{ [1]_{\rm res}\sqrt \vr\Big|\partial_t \vc U + \vc U \cdot \Grad \vc U\Big| \sqrt \vr \Big|\vc V - \vc v \Big|}{\rm d}t
\nonumber\\
&\qquad\le
\int_0^\tau\Big\|\partial_t\vc U+\vc U\cdot\Grad\vc U\|_{L^\infty(\Omega;\R^3)}\Big\|\Big[\vr\Big]_{\rm res}\|_{L^1(\Omega)}^{1/2}
\Big\|\vr \Big({\vc v}-\vc V\Big)^2\Big\|_{L^1(\Omega)}^{1/2} {\rm d}t
\nonumber\\
&\qquad\le c(\underline r,\overline r)\int_0^\tau a(t) {\cal E}\Big(\vr,\vc v\Big| r,\vc V\Big){\rm d}t
\end{align*}
with the same $a$ as before. In all above three formulas, we have employed (\ref{rentropy}) in the passage to  their last lines.

As far as the third term is concerned, we have immediately,
\begin{equation}\label{remainder3}
\int_0^\tau\intO{
\vr(\vc v- \vc V) \cdot \Grad \vc U
\cdot (\vc V - \vc v)}{\rm d}t \le c  \int_0^\tau a(t)
{\cal E}\Big(\vr,\vc v\Big| r,\vc V\Big){\rm d}t
\end{equation}
with
$$
a=\|\nabla\vc U\|_{L^\infty(\Omega;R^{9})}\in L^2(0,T).
$$

Resuming, the first volume integral in the remainder ${\cal R}$ of (\ref{reis}) is bounded from above by
\bFormula{remainder1st}
 \delta\int_0^\tau\|\vc v-\vc V\|^2_{W^{1,2}(\Omega;\R^3)}{\rm d} t+ c  \int_0^\tau a(t)
{\cal E}\Big(\vr,\vc v\Big| r,\vc V\Big){\rm d}t,
\eF
where $\delta>0$, $c=c(\delta,\underline r,\overline r)>0$ and
$$
a= \|\nabla\vc U\|_{L^\infty(\Omega;\R^{9})} + \|\partial_t\vc U+\vc U\cdot\Grad\vc U\|_{L^\infty(\Omega;\R^3)}^2\in L^1(0,T).
$$
{\bf Step 3:}{\it The second volume integral in ${\cal R}$}\\
As far as the last term is concerned, we use: 1) The Taylor formula together with the regularity $C^2$ of the pressure $p$, in order to estimate the essential part
\begin{align}\label{remainder5}
&- \int_0^\tau\intO { \Big[p(\rho)-p'(r)(\rho-r) -p(r)\Big]_{\rm ess} {\rm div} \vc U }{\rm d} t\nonumber\\
&\qquad
\le
c(\underline r,\overline r, |p'|_{C^1([\underline r/2,2\overline r])})\int_0^\tau \|{\rm div}\vc U\|_{L^\infty(\Omega)}\Big\|\Big[\vr- r\Big]_{\rm ess}\Big\|^2_{L^2(\Omega)}
\nonumber\\
&\qquad\le
c(\underline r,\overline r, |p'|_{C^1([\underline r/2,2\overline r])}) \int_0^\tau a(t)
{\cal E}\Big(\vr,\vc  v\Big| r,\vc V\Big){\rm d}t,\quad a=\|{\rm div}\vc U\|_{L^\infty(\Omega)}\in L^2(0,T).
\end{align}
2) Employing Lemma \ref{Lrelaxed2}, hypotheses (\ref{w1+}) and Corollary \ref{Crelaxed3}, we deduce the pointwise bound,
$$
\Big|[p(\rho)-p'(r)(\rho-r) -p(r)]_{\rm res}\Big|\le c(\underline r,\overline r, |p|_{C^1([\underline r,\overline r])}) E(\vr\Big|r)
$$
in order to estimate the residual part
\begin{align}\label{remainder6}
&- \int_0^\tau\intO { \Big[p(\rho)-p'(r)(\rho-r) -p(r)\Big]_{\rm res} {\rm div} \vc U }{\rm d} t
\le
c (\underline r,\overline r, |p|_{C^1([\underline r,\overline r])}) \int_0^\tau a(t)
{\cal E}\Big(\vr,\vc v\Big| r,\vc V\Big){\rm d}t.
\end{align}

We resume estimates obtained in Step 2:
\begin{equation}\label{remainder2nd}
- \int_0^\tau\intO { \Big[p(\rho)-p'(r)(\rho-r) -p(r)\Big]_{\rm ess} {\rm div} \vc U }{\rm d} t\le
c\int_0^\tau a(t){\cal E}\Big(\vr,\vu\Big| r,\vc V\Big){\rm d} t,
\end{equation}
where
$$
a=\|{\rm div}\vc U\|_{L^\infty(\Omega)}\in L^2(0,T), \; c= c(\underline r,\overline r, |p|_{C^2([\underline r/2,2\overline r])})>0.
$$
{\bf Step 4:}{\it The third volume integral in ${\cal R}$}\\
Similarly as in Step 1,
\begin{align}\label{remainder4}
& \int_0^\tau\intO{ \frac{\Grad p(r)}{r} ( r-\rho)\cdot( \vc v - \vc V)} {\rm d}t
\nonumber\\
&\qquad\le \delta \int_0^\tau \Big\|{\vc v}-\vc V\Big\|^2_{W^{1,2}(\Omega;\R^3)}{\rm d}t + c(\delta,\underline r, \overline r) \int_0^\tau a(t)
{\cal E}\Big(\vr,\vc v\Big| r,\vc V\Big){\rm d}t,
\end{align}
where
$$
a=\Big\| \frac{\Grad p(r)}{r}\Big\|^2_{L^\infty(\Omega;\R^3)}\in L^1(0,T).
$$
{\bf Step 5:}{\it Conclusion}\\
Coming back with these estimates to the relative energy inequality (\ref{reis}), taking into account (\ref{belowbound}) and choosing $\delta$ sufficiently small with respect to $\mu$, we easily verify the validity of (\ref{how3}). This finishes the proof of Theorem \ref{TWS1}.

\section{Concluding remarks}\label{CR}
\subsection{Existence of dissipative solutions in domains with piecewise regular boundaries}
So far we have established the existence of dissipative weak solutions  under the
assumption of the 𝐶$C^2$-regularity of the domain. In many practical situations in nonzero outflow/inflow regimes, the domain
occupied by the fluid does not possess this regularity. A typical example of such situation is a finite cylinder with inflow and
outflow boundaries lower and upper discs of the boundary of the cylinder. The present section intends to remove this drawback.

We start with the definition of the piecewise $C^2$ Lipschitz domain.
\\
{\bf Definition \ref{CR}.1}[Piecewise $C^2$ Lipschitz domain] \\
{\it
We shall say that $\Omega\subset R^d$, $d=2,3$, is a bounded piecewise $C^{2}$ Lipschitz domain if
\begin{description}
\item{\it 1.} $\Omega$ is a bounded  Lipschitz domain;
\item{\it 2.}
The boundary of the domain can be written as
\[
\partial\Omega=\Gamma\cup\gamma \quad \text{ with }
\Gamma=\cup_{i=1}^{I}\Gamma_i, \quad
\gamma=\cup_{i=1}^{I} \gamma_i
\]
where $\Gamma_i$ are open connected $(d-1)$ dimensional mutually disjoint manifolds of class $C^{2}$ and
\[
\gamma_i\equiv\partial\Gamma_i=\cup_{k_i=1}^{K_{i}}\gamma_{i,k_i}.
\]
When $d=3$, $\gamma_{i,k_i}$ is a closed parametrized curve in $R^3$ of class $C^{2}$.
{ If $\gamma_{i,k_i}$ and $\gamma_{j,l_j}$ intersect, they coincide.}
When $d=2$, $\gamma_{i,k_i}$ is a point in $R^2$.
\end{description}
}

\bTheorem{M2}
Let $\vr_B$, $\vu_B$, $\vr_0$, $\vu_0$ and $p$ satisfy all assumptions of Theorem \ref{TM1}.
We suppose that $\Omega$ is a piecewise $C^2$ Lipschitz domain such that:
\begin{description}
\item{\it 1.}
\bFormula{dom1}
\partial\Omega
=\overline\Gamma_0\cup\overline\Gamma_{\rm in}\cup\overline \Gamma_{\rm out}.
\eF
\item{\it 2.}
There holds
\bFormula{dom2}
\begin{split}
&\Gamma_{\mathfrak a}= \cup_{i_{\mathfrak a}=1}^{I_{\mathfrak a}}\Gamma_{i_{\mathfrak a}} \quad \text{ where ${\mathfrak a}$ stands for { ``$0$'', {\rm ``in''}, {\rm ``out''}}}, \\
&\overline\Gamma_{k_{\mathfrak{ a}}}\cap \overline\Gamma_{l_{\mathfrak{b}}}=\emptyset \quad \text{ whenever $\mathfrak{a}\in\{{\rm in},{\rm out}\}$,  $\mathfrak{ b}\in\{{\rm in},{\rm out}\}$, $ k_{\mathfrak{ a}}\neq l_{\mathfrak{ b}}$},
\end{split}
\eF
where $\Gamma_{i_{\mathfrak a}}$ are  (open, connected) $(d-1)$-dimensional mutually disjoint manifolds of class $C^{2}$,  { ``in'' and ``out'' refer} to the notation \eqref{i5} and
\bFormula{G0}
\Gamma_0={\rm int}_{d-1}\Big(\{x\in \partial\Omega\,|\,\vu_B\cdot\vc n=0\}\Big).
\eF
In the above ${\rm int}_{d-1}$ the interior on the (hyper)surface $\partial\Omega$.
\item{\it 3.} There holds
\bFormula{dom3}
\gamma_{\mathfrak a}\equiv\partial\Gamma_{\mathfrak a}=\cup_{k_{\mathfrak a}=1}^{K_{\mathfrak a}}\gamma_{{\mathfrak a},k_{\mathfrak a}},
\eF
where $\gamma_{{\mathfrak a},k_{\mathfrak a}}$ is a closed parametrized  curve in $R^d$ of class $C^{2}$ (if $d=3$) or a point (if $d=2$)
such that
either $\gamma_{{\mathfrak a}, k_{\mathfrak a}}\cap\gamma_{{\mathfrak b},l_{\mathfrak b}}= \emptyset$ or $\gamma_{{\mathfrak a},k_{\mathfrak a}}=\gamma_{{\mathfrak b},l_{\mathfrak b}}$.
\end{description}
Then all conclusions of Theorem \ref{TM1} remain valid, in particular, the problem (\ref{i1}--\ref{i5}) admits a dissipative weak solution.
\eT
\bRemark{piecewise}
\begin{enumerate}
\item
The main issue of the proof of Theorem \ref{TM2} is a construction of a convenient approximation of the piecewise regular domain $\Omega$ by regular ($C^2$) bounded domains $\Omega_\kappa$ (small parameter $\kappa>0$) keeping conserved up to { small perturbations}  the inflow/outflow properties of the fluid flow. Such approximation of the domain and boundary data has been suggested in \cite[Section 3]{ChoeNoY} and
existence of bounded energy weak solutions has been proved through the limit passage $\kappa\to 0$
from bounded energy weak solutions on $C^2$-domains $\Omega_\kappa$ to bounded energy weak solutions on domain $\Omega$. Likewise, using ideas
of Section \ref{D}, one can pass to the limit also from the relative energy inequality on $\Omega_\kappa$ to relative energy inequality
on $\Omega$. We let the details to the interested reader.
\item The dissipative weak solutions constructed in Theorem \ref{TM2} fulfill all assumptions of Theorem \ref{TWS1}. In particular, they satisfy the weak-strong uniqueness principle.
\end{enumerate}
\eR
\subsection{Nonmonotone pressure law}

The statement of Theorem \ref{TM1} and also of Theorem \ref{TM2} can be generalized to some  possibly non monotone pressure laws, as, e.g.,
\begin{equation}\label{pressure1+}
p=\pi+\mathfrak{p},\;
\pi\in C[0,\infty)\cap C^1(0,\infty), \;\pi(0)=0, \; {\pi}'(\vr)>0,
\end{equation}
$$
\pi'(\vr)\ge a_1\vr^{\gamma-1}-b,\; \pi(\vr)\le a_2\vr^\gamma+ b
$$
with $\gamma>d/2$ and $a_1,a_2>0$, $b\ge 0$, and
\begin{equation}\label{pp}
\mathfrak{p}\in C^2_c[0,\infty),\; \mathfrak{p}\le 0.
\end{equation}
In this case the existence of bounded energy weak solutions in non zero inflow/outflow setting has been proved in \cite{JiNoSIMA} and \cite{ChoeNoY}
and consequently the construction of dissipative weak solution can be performed by pursuing the strategy of the present paper.

Likewise the statement of Theorem \ref{TWS1} can be generalized to the pressure laws (\ref{pressure1+}) with $\gamma>1$ and with
\begin{equation}\label{pp+}
\mathfrak{p}\in C^2_c[0,\infty)\,\mbox{globally Lipschitz}.
\end{equation}
In this case the relative energy function $E(\vr|r)$ defined in (\ref{E}) must be calculated from the monotone part of the pressure, i.e., with
$$
H(\vr)=\vr\int_1^\vr\frac{\pi(z)}{z^2}{\rm d} z.
$$
This result can be proved by combining  the strategy of the present paper with the ideas introduced in Feireisl \cite{Fnmws} and
Chaudhuri \cite{Chaud}.

\def\cprime{$'$} \def\ocirc#1{\ifmmode\setbox0=\hbox{$#1$}\dimen0=\ht0
  \advance\dimen0 by1pt\rlap{\hbox to\wd0{\hss\raise\dimen0
  \hbox{\hskip.2em$\scriptscriptstyle\circ$}\hss}}#1\else {\accent"17 #1}\fi}


\begin{thebibliography}{10}


\bibitem{BreJab}
{\sc D.~Bresch and P.-E. Jabin,}
\newblock{\em Global existence of weak solutions for compresssible
  {N}avier-{S}tokes equations: {T}hermodynamically unstable pressure and
  anisotropic viscous stress tensor,}
\newblock{Ann. of Math.}, { 188} (2018), pp.-577--684.

\bibitem{BG10}
{\sc S. Benzoni-Gavage,}
\newblock{\em Calcul diff\'erentiel et \'equations diff\'erentielles,}
\newblock{Dunod,} 2010.




\bibitem{JiNoSIMA}
{\sc T. Chang, B.J. Jin, A. Novotny}
\newblock
{\it Compressible Navier-Stokes system with general inflow-outflow
boundary data}
\newblock{ SIAM J. Math. Anal., 2019}, accepted

\bibitem{Chaud}
{\sc N. Chaudhuri}
\newblock {\em On weak-strong uniqueness for compressible Nvier-Stokes system with general pressure laws}
\newblock{NORWA, 2019}, accepted

\bibitem{ChoeNoY}
{\sc H.J. Choe, A. Novotny, M. Yang}
\newblock		{\it Compressible Navier-Stokes system with general inflow-outflow
boundary data on piecewise regular domains}
\newblock{ ZAMM Z. Angew. Math. Mech. 98(8) (2018), pp. 1447--1471}

\bibitem{Dafermos}
{\sc C.M. Dafermos.}
\newblock{\em The second law of thermodynamics and stability.}
\newblock Arch. Rational Mech. Anal. 70 (1979), pp. 167-179


\bibitem{DHP}
{\sc R. Denk, M. Hieber, and J.~Pr{\" u}ss,}
\newblock{\em Optimal $Lp-Lq$-estimates for parabolic boundary value
problems with inhomogeneous data},
\newblock{ Math. Z.}, {257}(2007), pp.-193–-224.



\bibitem{DL}
{\sc R.J. DiPerna and P.-L. Lions},
\newblock{\em Ordinary differential equations, transport theory and {S}obolev
  spaces},
\newblock { Invent. Math.}, {98}(1989), pp.-511--547.

\bibitem{Evans}
{\sc L. C. Evans},
\newblock {\em Partial Differential Equations},
\newblock Graduate Studies in Mathematics, Vol. 19, AMS

\bibitem{EF70}
{\sc E.~Feireisl},
\newblock {\em Dynamics of viscous compressible fluids},
\newblock Oxford University Press, Oxford, 2004.

\bibitem{FEnonmon}
{\sc  E. Feireisl},
  \newblock{\em Compressible Navier-Stokes equations with a non-monotone pressure law,}
  \newblock{J. Differential Equations}, {184} (2002), pp.-97--108.
	
	\bibitem{Fnmws}
	{\sc E. Feireisl},
	\newblock{\em On weak-strong uniqueness for the compressible Navier-Stokessystem with non-monotone pressure law},
	\newblock{Preprint 27-2018}, Math. Inst. Czech. Acad. Sci, \url{http://www.math.cas.cz/fichier/preprints/IM_20180623094424_14.pdf}
	
	\bibitem{FeJiNo}
{\sc E.~Feireisl, B~J.~Jin, and A.~Novotn{\'y}}
\newblock Relative entropies, suitable weak solutions, and weak-strong
  uniqueness for the compressible Navier--Stokes system.
\newblock {\em J. Math. Fluid. Mech.}, 14(4):717--730, 2012.

\bibitem{FEINOV}
{\sc E.~Feireisl and A.~Novotn\'y,}
\newblock {\em Singular limits in thermodynamics of viscous fluids},
\newblock Birkhauser, Basel, 2009. 2nd. enlarged edition 2017



\bibitem{FNP}
{\sc E.~Feireisl, A.~Novotn{\' y} and H.~Petzeltov{\' a}},
\newblock{\em On the existence of globally defined weak solutions to the
  {N}avier-{S}tokes equations of compressible isentropic fluids},
\newblock { J. Math. Fluid Mech.}, { 3}(2001), 358--392.

\bibitem{FNSun}
{\sc E. Feireisl, A. Novotn\'y, Y. Sun}
\newblock{Suitable weak solutions to the Navier-Stokes equations of compressible viscous fluids}
\newblock { Indiana Univ. Math. J. 60,  2011, 611 - 632} 


\bibitem{GaHeMaNo}
{\sc T. Gallouet, R. Herbin, D. Maltese, A. Novotny}
\newblock {\it Error estimates for a numerical approximation to the compressible barotropic Navier-Stokes equations},
  \newblock  { IMA Journal of Numerical Analysis, 36(2) (2016), pp. 543-592}
	
	\bibitem{GaMaNo}
	{\sc T. Gallouet, D. Maltese, A. Novotny}
	\newblock{\it  Error estimates for the implicit MAC scheme for the compressible Navier–Stokes equations.},
	\newblock{  Numer. Math. 141(2), 2019, 495–-567}
	
	\bibitem{Kleinetal}
	{\sc R. Klein, N. Botta, T. Schneider, C.D. Munz, S. Roller, A. Meister, L. Hoffman, T. Sonar}
	\newblock{\em Asymptotic adaptive methods for multi-scle problems in fluid mechanics.}
	\newblock{J. Eng. Math. 39 (2001), pp. 537-559}
	
	\bibitem{Galdi}
	{\sc G. P. Galdi}
	\newblock{An introduction to the mathematical theory of the Navier--Stokes equations I.}
	\newblock Springer, New-York, Second edition, 2011
	
	
\bibitem{Girinon}
{\sc V. Girinon},
\newblock{\em Navier-Stokes equations with nonhomogeneous boundary conditions in a bounded three-dimensional domain,}
\newblock{J. Math. Fluid Mech. }, 13(2011), pp.-309–-339.

\bibitem{LI4}
{\sc P.-L. Lions},
\newblock {\em Mathematical topics in fluid dynamics, Vol.2, Compressible
  models},
\newblock Oxford Science Publication, Oxford, 1998.




\bibitem{NOVO}
{\sc S.~Novo,}
\newblock{\em Compressible {N}avier-{S}tokes model with inflow-outflow boundary
  conditions},
\newblock {J. Math. Fluid Mech.}, {7}(2005), pp.-485--514.


\bibitem{NOST}
{\sc A.~Novotn{\' y} and I.~Stra{\v s}kraba},
\newblock {\em Introduction to the mathematical theory of compressible flow},
\newblock Oxford University Press, Oxford, 2004.









\bibitem{VAZA}
{\sc A.~Valli and M.~Zajaczkowski,}
\newblock{\em {N}avier-{S}tokes equations for compressible fluids: Global existence
  and qualitative properties of the solutions in the general case,}
\newblock { Comm. Math. Phys.}, {103}(1986), pp.-259--296, 1986.

\end{thebibliography}
\end{document}